\documentclass[reqno, 11pt, a4paper]{amsart}
\usepackage{latexsym,ifthen,xspace}
\usepackage{amsmath,amssymb,amsthm}
\usepackage{dsfont}
\usepackage{mathtools}
\usepackage{enumerate, calc, bbm}
\usepackage[usenames,dvipsnames]{xcolor}
\usepackage[colorlinks=true, pdfstartview=FitV, linkcolor=blue, citecolor=blue, urlcolor=blue, pagebackref=false]{hyperref}
\usepackage{orcidlink}
\usepackage[text={33pc,605pt},centering]{geometry} 
\usepackage{graphicx, xcolor}
\usepackage{array}
\usepackage{url}

\numberwithin{equation}{section}

\newtheorem{theorem}{Theorem}[section]
\newtheorem{lemma}[theorem]{Lemma}
\newtheorem{proposition}[theorem]{Proposition}
\newtheorem{corollary}[theorem]{Corollary}

\theoremstyle{definition}
\newtheorem{definition}[theorem]{Definition}
\newtheorem{assumption}[theorem]{Assumption}
\newtheorem{example}[theorem]{Example}
\newtheorem{notation}{Notation}

\theoremstyle{remark}
\newtheorem{remark}[theorem]{Remark}

\DeclareMathAlphabet{\mathsl}{OT1}{cmss}{m}{sl}
\SetMathAlphabet{\mathsl}{bold}{OT1}{cmss}{bx}{sl}

\DeclarePairedDelimiter{\norm}{\lVert}{\rVert}
\DeclarePairedDelimiter{\scpr}{\langle}{\rangle}
\DeclareMathOperator{\capa}{cap}

\DeclareMathOperator{\mean}{\mathbb{E}}
\DeclareMathOperator{\Mean}{\mathrm{E}}
\newcommand{\tMean}{\ensuremath{\tilde{\Mean}\mspace{-12mu}\phantom{\Mean}}}
\DeclareMathOperator{\prob}{\mathbb{P}} 
\DeclareMathOperator{\Prob}{\mathrm{P}}
\DeclareMathOperator{\Var}{\mathbb{V}ar}

\newcommand{\tProb}{\ensuremath{\widetilde{\Prob}\mspace{-12.5mu}\phantom{\Prob}}}
\newcommand{\cA}{\ensuremath{\mathcal A}}
\newcommand{\cB}{\ensuremath{\mathcal B}}
\newcommand{\cG}{\ensuremath{\mathcal G}}
\newcommand{\cM}{\ensuremath{\mathcal M}}
\newcommand{\cS}{\ensuremath{\mathcal S}}
\newcommand{\cX}{\ensuremath{\mathcal X}}
\newcommand{\me}{\ensuremath{\mathrm{e}}}
\newcommand{\md}{\ensuremath{\mathrm{d}}}
\newcommand{\ind}{\mathbbm{1}}
\newcommand{\abs}[1]{\left\lvert #1 \right\rvert}
\newcommand{\newquote}[1]{``#1''}
\newcommand{\Ometa}{\Omega_{\mathrm{meta}}}
\newcommand{\Otimeta}{\widetilde{\Omega}_{\mathrm{meta}}}

\begin{document}

\title[]{Metastability of Glauber dynamics with\\ 
  inhomogeneous coupling disorder}

\author[A. Bovier]{Anton Bovier}
\address{A. Bovier - University of Bonn\\ Institute for Applied Mathematics, Endenicher Allee 60, 53115 Bonn, Germany}
\email{bovier@uni-bonn.de}

\author[F. den Hollander]{Frank den Hollander}
\address{F. den Hollander - Leiden University\\
  Mathematical Institute\\ Einsteinweg 55, 2333CC Leiden, The Netherlands}
  \email{denholla@math.leidenuniv.nl}

\author[S. Marello]{Saeda Marello}
\address{S. Marello - University of Bonn\\ Institute for Applied Mathematics \\ Endenicher Allee 60, 53115 Bonn, Germany}
\email{marello@iam.uni-bonn.de}

\author[E. Pulvirenti]{Elena Pulvirenti}
\address{E. Pulvirenti - Delft University of Technology\\
  Delft Institute of Applied Mathematics (DIAM)\\ Mekelweg 4, 2628CD Delft, Netherlands}
\email{e.pulvirenti@tudelft.nl}

\author[M. Slowik]{Martin Slowik}
\address{M. Slowik - University of Mannheim\\
  Mathematical Institute\\ B6, 26, 68159 Mannheim, Germany}
\email{slowik@math.uni-mannheim.de}

\begin{abstract}

  We introduce a general class of mean-field-like spin systems with random couplings that comprises both the Ising model on inhomogeneous dense random graphs and the randomly diluted Hopfield model.  We are interested in quantitative estimates of metastability in large volumes at fixed temperatures when these systems evolve according to a Glauber dynamics, i.e.\ where spins flip with Metropolis transition probabilities at inverse temperature $\beta$.  We identify conditions ensuring that with high probability the system behaves like the corresponding system where the random couplings are replaced by their averages.  More precisely, we prove that the metastability of the former system is implied with high probability by the metastability of the latter.  Moreover, we consider relevant metastable hitting times of the two systems and find the asymptotic tail behaviour and the moments of their ratio.  This work provides an extension of the results known for the Ising model on the Erd\H{o}s--R{\'e}nyi random graph.  The proofs use the potential-theoretic approach to metastability in combination with concentration inequalities.

\end{abstract}

\thanks{AB and SM were partly funded by the Deutsche Forschungsgemeinschaft (DFG) under Germany's Excellence Strategy - GZ 2047/1, Projekt-ID 390685813, and through Project-ID 211504053 - SFB 1060. FdH was partly funded by the Netherlands Organisation for Scientific Research (NWO) through Gravitation-grant NETWORKS-024.002.003.}

\date{\today}

\maketitle


\section{Introduction}


Over the last decade there has been increasing interest in metastability under Glauber dynamics of the Ising model with \emph{random interactions}, in particular, of the Ising model on random graphs. Dommers \cite{Do17} considered the case of random regular graphs, Dommers, den Hollander, Jovanovski, and Nardi \cite{DdHJN17} the configuration model, in both cases in finite volume and at low temperature.  Mossel and Sly \cite{MS09, MS13} computed mixing times on sparse Erd\H{o}s-R{\'e}nyi random graphs and on random regular graphs, in both cases in finite volume and at high temperature.  Recently, Can, van der Hofstad, and Kumagai \cite{CvdHK21} analysed mixing times on random regular graphs, in large volumes and at fixed temperature. 

Metastability under Glauber dynamics of the Ising model on dense random graphs has so far only been studied for the Erd\H{o}s-R{\'e}nyi random graph with fixed edge retention probability, by den Hollander and Jovanovski \cite{dHJ21} and by Bovier, Marello, and Pulvirenti \cite{BMP21}.  In both papers, mean metastable exit times of the random model are compared to those of the standard Curie-Weiss model, in large volumes and at fixed temperature.  In \cite{dHJ21} the pathwise approach to metastability (see Olivieri and Vares \cite{OV05}) was used to prove that mean metastable exit times are asymptotically equal to those of the Curie-Weiss model, multiplied by a random prefactor of polynomial order in the size of the system. The prefactor estimate was improved in \cite{BMP21} by using the potential-theoretic approach to metastability (see Bovier and den Hollander \cite{BdH15}), at the expense of losing generality in the initial distribution. Recently, Bovier, den Hollander and Marello \cite{BdHM22} studied metastability under Glauber dynamics of the Ising model on the complete graph with random independent couplings in large volumes and at fixed temperature. In that model, the product structure of the couplings allows for lumping of states, and for combining the potential-theoretic approach with coarse-graining techniques, to obtain sharp estimates on mean metastable exit times.

The present paper extends the  results for the Erd\H{o}s--R{\'e}nyi random graph to inhomogeneous dense random graphs and to more general random interactions. We compare the metastable behaviour of a class of spin systems whose Hamiltonian has random and conditionally independent coupling coefficients, called \emph{quenched model}, with the corresponding \emph{annealed model} in which the coupling coefficients are replaced by their conditional mean. More precisely, we prove that metastability of the annealed model implies, in large volumes and at fixed temperature, almost sure metastability of the quenched model with respect to the metastable sets of the annealed model. Moreover, assuming metastability of the annealed model, we consider the ratio between the mean hitting times of the quenched model and the annealed model, and estimate both its tail behaviour and its moments, again in large volumes and at fixed temperature.

As in \cite{BMP21}, we follow the potential-theoretic approach to metastability, which\break allows us to estimate mean metastable exit times by estimating \emph{capacities} and\break weighted sums of the equilibrium potential called \emph{harmonic sums}. Estimates on the former can be obtained with the help of well-known variational principles, while estimates on the latter are more involved. See, for instance, Bianchi, Bovier and Ioffe \cite{BBI09} and \cite{BMP21}, where long and model-dependent computations were needed to prove that the relevant contribution of the harmonic sum is localised around the starting metastable set.  Schlichting and Slowik in \cite{SS19}, using an  alternative definition of metastable sets,  prove that localisation of the harmonic sum around the starting metastable set holds in large generality. Their work allows us to derive  results for a large class of models. A second novelty of the present paper compared to \cite{BMP21} concerns the techniques that are used to prove concentration results. In \cite{BMP21}, Talagrand's concentration inequality was used, while here we use McDiarmid's concentration inequality.

\section{Model, results and methods}
This section is structured as follows. In Section~\ref{sec:model}, we introduce the model. In Section~\ref{sec:meta}, we define metastability, introduce relevant quantities, and state our main results. In Section~\ref{sec:meth}, we summarise our strategy and methods, and give an outline of the rest of the paper.


\subsection{The model}
\label{sec:model}

Let $(\Omega, \mathcal{F}, \prob)$ be an abstract probability space. Let $\mathcal{G} \subset \mathcal{F}$ be a sub-$\sigma$-algebra of $\mathcal{F}$ and $J = (J_{ij})_{1 \leq i < j < \infty}$ be a triangular array of real random variables that are \emph{conditionally independent} given $\mathcal{G}$ and \emph{uniformly bounded}, i.e., there exists a $k_J \in (0, \infty)$ such that $\vert J_{ij}\vert  \leq k_J$ $\prob$-a.s.\ for all $1 \leq i < j < \infty$.  We write $\prob_{\mathcal{G}}[\,\cdot\,]$ to denote a regular conditional distribution for $J$ given $\mathcal{G}$ (which exists because $J$ is a sequence of real random variables; see Chow and Teicher \cite[p.~218]{CT97}).  Write $\mean$ to denote the expectation with respect to $\prob$, and $\mean_{\mathcal{G}}$ and $\Var_{\mathcal{G}}$ to denote expectation and variance with respect to $\prob_{\mathcal{G}}$.

Given $2 \leq N \in \mathbb{N}$, consider the Ising Hamiltonian with random couplings of the form
\begin{align}
  \label{eq:def_HN}
  H_N(\sigma) 
  \coloneqq
  -\frac{1}{N} \sum_{\substack{i,j = 1 \\ i<j}}^N J_{ij} \sigma_i \sigma_j
  - h \sum_{i=1}^N \sigma_i,
  \qquad \sigma \in \cS_N,
\end{align} 
with $h \in \mathbb{R}$ the magnetic field and $\cS_N\coloneqq \{-1,1\}^N$ the set of spin configurations. The corresponding Gibbs measure on $\cS_N$ is denoted by 
\begin{align}
  \label{eq:def_muN}
  \mu_N(\sigma)
  \coloneqq 
  Z_N^{-1}\, \me^{-\beta H_N(\sigma)},
  \qquad \sigma \in \cS_N,
\end{align}
with $\beta \in (0, \infty)$ the inverse temperature and $Z_N$  the normalizing partition function. The spin configurations evolve in time as a discrete-time Markov chain $(\Sigma_N(t))_{t \in \mathbb{N}_0}$ with state space $\cS_N$ and Glauber-Metropolis transition probabilities given by 
\begin{align}
  \label{eq:defpN}
  p_N(\sigma, \sigma')
  \coloneqq 
  \begin{cases}
    \frac{1}{N}\exp\bigl(-\beta \bigl[H_N(\sigma') - H_N(\sigma)\bigr]_+\bigr),
    \quad 
    &\text{if } \sigma \sim \sigma', \\
    1 - \sum_{\eta \neq \sigma} p_N(\sigma, \eta),
    &\text{if } \sigma = \sigma', \\
    0, 
    &\text{otherwise},
  \end{cases}
  \qquad \sigma, \sigma' \in \cS_N,
\end{align}
where $\sigma \sim \sigma'$ means that $\sigma'$ is obtained from $\sigma$ by a flip of a single spin. The associated (discrete) generator $L_N$ acts on bounded functions $f :\cS_N \to \mathbb{R}$ as 
\begin{align}
  \bigl(L_N f\bigr)(\sigma)
  \coloneqq 
  \sum_{\sigma' \in \cS_N} p_N(\sigma, \sigma') 
  \bigl(f(\sigma') - f(\sigma)\bigr), \qquad \sigma \in \cS_N.
\end{align}

Note that the stochastic process $(\Sigma_N(t))_{t \in \mathbb{N}_0}$ is irreducible and reversible with respect to the Gibbs measure $\mu_N$.  We denote by $\Prob_{\nu}^N$ the law of the Markov chain $(\Sigma_N(t))_{t \in \mathbb{N}_0}$ with initial distribution $\nu$. The corresponding expectation is denoted by $\Mean_{\nu}^N$. If the initial distribution is concentrated on a single configuration $\sigma \in \cS_N$, then we write $\Prob_{\sigma}^N$ and $\Mean_{\sigma}^N$, respectively. For a non-empty subset $\mathcal{A} \subset \cS_N$, let $\tau_{\mathcal{A}}^N$ be the first return time to $\mathcal{A}$, i.e., 
\begin{align}
  \tau_{\mathcal{A}}^N
  \equiv
  \tau_{\mathcal{A}}^N\bigl((\Sigma_N(t))_{t \in \mathbb{N}_0}\bigr)
  \coloneqq 
  \inf\bigl\{
    t \in \mathbb{N}: \Sigma_N(t) \in \mathcal{A}
  \bigr\}.
\end{align}

Our main objective is to compare the evolution of the model with Hamiltonian $H_N$ and the model with Hamiltonian $\tilde{H}_N$ defined by
\begin{align}
  \label{eq:def_EHN}
  \widetilde{H}_N(\sigma)
  \coloneqq 
  \mean_{\mathcal{G}} \bigl[H_N(\sigma)\bigr]
  =
  -\frac{1}{N} \sum_{\substack{i,j = 1 \\ i<j}}^N 
  \mean_{\mathcal{G}}\bigl[J_{ij}\bigr] \sigma_i \sigma_j
  - h \sum_{i=1}^N \sigma_i,
  \quad \sigma \in \cS_N,
  \quad \prob\text{-a.s.}
\end{align}
Throughout the paper, we use the superscript $\sim$ to denote quantities that refer to the model defined in terms of $\widetilde{H}_N$. For instance, $((\widetilde{\Sigma}_N(t))_{t \in \mathbb{N}_0}, \tProb_{\sigma}^N\colon \sigma \in \cS_N)$ denotes the discrete-time Markov chain with transitions probabilities $(\tilde{p}_N(\sigma, \sigma'))_{\sigma, \sigma' \in \cS_N}$, which is reversible with respect to $\tilde{\mu}_N$, where both the transition probabilities and the Gibbs measure are defined as in \eqref{eq:def_HN} and \eqref{eq:def_muN}, but with $H_N$ replaced by $\widetilde{H}_N$.  For lack of better names and with abuse of terminology, we  refer to the models defined in terms of $H_N$ and $\widetilde{H}_N$ as the \emph{quenched model} and the \emph{annealed model}, respectively.

In the sequel, we provide a list of motivating examples for which the results stated later hold. For this purpose, consider two sequences 
\begin{equation}
  (A_{ij})_{1 \leq i < j < \infty}, \qquad (P_{ij})_{1 \leq i < j < \infty},
\end{equation}
of triangular arrays with $\vert A_{ij}\vert  \leq k_J$ and $P_{ij} \in (0,1)$ for $i,j \in \mathbb{N}$ with $i < j$, and let $\mathcal{G} \coloneqq \sigma(A_{ij}, P_{ij} \colon 1 \leq i < j \leq \infty)$ be the $\sigma$-algebra generated by these sequences.  Moreover, let 
\begin{equation}
  (U_{ij})_{1 \leq i < j < \infty}
\end{equation}
be a triangular array of i.i.d.\ random variables distributed uniformly in $(0,1)$. Define 
\begin{align}
  J_{ij} 
  \coloneqq 
  A_{ij} B_{ij}, 
  \qquad
  B_{ij} 
  \coloneqq 
  \ind_{\{U_{ij} \leq P_{ij}\}},
  \qquad 1 \leq i < j <\infty.
\end{align}
Note that $(J_{ij})_{1 \leq i < j < \infty}$ and $(B_{ij})_{1 \leq i < j < \infty}$ are triangular arrays of \emph{conditionally independent} random variables given $\mathcal{G}$. In particular, $B_{ij}$ are Bernoulli random variables with mean $P_{ij}$.

\begin{example}[Ising model on the Erd\H{o}s--R{\'e}nyi random graph]
  \label{ex:dilute:CW}
  By choosing $A_{ij} \coloneqq 1$ and $P_{ij} \coloneqq p \in (0,1]$ for $1 \leq i < j < \infty$, $\mathcal{G}$ becomes the set $\{\emptyset,\Omega\}$ and $H_N$ in \eqref{eq:def_HN} becomes the Hamiltonian of the Ising model on the Erd\H{o}s--R{\'e}nyi random graph with edge retention probability $p$, known as the \emph{randomly diluted Curie-Weiss model}. Its metastable behaviour was studied in \cite{dHJ21} and \cite{BMP21}. In this case the annealed model is the Curie-Weiss model.
\end{example}  

\begin{example}[Ising model on inhomogeneous random graphs]
  By taking $A_{ij} \coloneqq 1$ for  $1 \leq i < j < \infty$, $H_N$ in \eqref{eq:def_HN} becomes the Hamiltonian of the Ising model on an inhomogeneous random graph, in which an edge $(ij)$ is present with probability $P_{ij}$. Of particular interest is the case $P_{ij} = V_i V_j$ with $(V_i)_{i \in \mathbb{N}}$ a sequence of i.i.d.\ random variables with support in $(0,1)$, known as the Ising model on the \emph{Chung-Lu random graph} \cite{CL02}. The metastable behaviour of the corresponding annealed model was studied in \cite{BdHM22} for the case where the random variables $V_i$ have finite support.
\end{example}

\begin{example}[Randomly diluted Hopfield model]
  Given $n \in \mathbb{N}$ random patterns $\xi^1, \dots, \xi^n$, with  $\xi^k=(\xi_i^k)_{i \in \mathbb{N}}$ and $\xi_i^k \in {[-1,1]}$ for $1\leq k \leq n$, set $A_{ij} \coloneqq \sum_{k=1}^n \xi^k_i \xi^k_j$.  By taking $P_{ij}\equiv p \in (0,1)$ for $1 \leq i < j < \infty$, $H_N$ in \eqref{eq:def_HN} becomes the Hamiltonian of a \emph{Hopfield model} in which the interaction coefficients are \emph{randomly diluted} by i.i.d.\ Bernoulli random variables with mean $p$. See Bovier and Gayrard \cite{BG98} for a review on the Hopfield model.  The metastable behaviour of the annealed model, i.e., the undiluted Hopfield model, was studied by an der Heiden in \cite{adH07} in a restricted $(\beta,h)$-regime. We plan to address the metastable behaviour in a more general $(\beta,h)$-regime in a future paper.
\end{example}


\subsection{Metastability and main results}
\label{sec:meta}

Before stating our main results, we recall the definition of metastable Markov chains and metastable sets put forward in Schlichting and Slowik \cite[Definition 1.1]{SS19}.

\begin{definition}[$\rho_N$-metastability and metastable sets]
  \label{def:meta}
  For $\rho_N>0$ and $2 \leq K \in \mathbb{N}$, let $\{\mathcal{M}_{1,N}, \dots, \mathcal{M}_{K,N}\}$ be a set of subsets of $\cS_N$ such that $\mathcal{M}_{i,N} \cap \mathcal{M}_{j,N} = \emptyset$ for all $1 \leq i \neq j \leq K$.  The Markov process $(\Sigma_N(t))_{t \in \mathbb{N}_0}$ is called \emph{$\rho_N$-metastable} with respect to $\{\mathcal{M}_{1,N}, \ldots, \mathcal{M}_{K,N}\}$ when
  \begin{align}
    \label{eq:def:meta}
    K\,
    \frac{
      \max_{j \in \{1, \ldots, K\}}
      \Prob_{\mu_N \vert \mathcal{M}_{j,N}}^N\Bigl[
        \tau_{\mathcal{M}_N \setminus \mathcal{M}_{j,N}}^N 
        < \tau_{\mathcal{M}_{j,N}}^N
      \Bigr]
    }{
      \min_{\cX \subset \cS_N \setminus \mathcal{M}_N }
      \Prob_{\mu_N\vert \cX}^N\Bigl[ 
        \tau_{\mathcal{M}_N}^N < \tau_{\cX}^N 
      \Bigr]
    }
    \leq
    \rho_N \ll 1,
  \end{align}
  where $\mathcal{M}_N \coloneqq \bigcup_{i=1}^K \mathcal{M}_{i,N}$ and, for a non-empty set $\cX \subseteq \cS_N$, $\mu_N \vert \cX$ denotes the Gibbs measure $\mu_N$ conditioned on the set $\cX$.
\end{definition}

\begin{remark}
  The advantage of this definition compared to the one given in Bovier and den Hollander \cite[Definition~8.2]{BdH15} is twofold: it allows for direct control of $\ell^p(\mu_N)$-norms of functions, and does not depend on the cardinality of the state space.  For a more detailed comparison of the two definitions of metastability we refer to \cite[Remark~1.2]{SS19}.
\end{remark}

For fixed $2 \leq K \in \mathbb{N}$ and $k_1 > 0$, define
\begin{align}
  \Otimeta(N)
  \coloneqq 
  \Bigl\{ 
    \omega \colon 
    \exists\, \{\mathcal{M}_{1,N}, \ldots, \mathcal{M}_{K,N}\}(\omega) \text{ non-empty disjoint subsets of } \cS_N 	
    \nonumber\\
    \text{ s.t. } 
    (\widetilde{\Sigma}_N(t))_{t \in \mathbb{N}_0}(\omega) \text{ is $\me^{-k_1 N}$-metastable} \text{ w.r.t. } \{\mathcal{M}_{1,N}, \ldots, \mathcal{M}_{K,N}\}(\omega)
  \Bigr\},
  \label{eq:def_Otimeta}
\end{align}
i.e., the event that the Markov chain $(\widetilde{\Sigma}_N(t))_{t \in \mathbb{N}_0}$ is $\tilde{\rho}_N$-metastable with respect to some $\{\mathcal{M}_{1,N}, \ldots, \mathcal{M}_{K,N}\}$, where we abbreviate
\begin{align}\label{eq:def_tirho}
  \tilde{\rho}_N\coloneqq \me^{-k_1 N}.
\end{align}

\begin{remark}
  \label{rem:Gmeas}
  Note that both $\Otimeta(N)$ and $\{\mathcal{M}_{1,N}, \ldots, \mathcal{M}_{K,N}\}$ -- playing the role of \emph{metastable sets} -- are $\cG$-measurable, because they are defined in terms of the annealed Hamiltonian $\widetilde{H}_N$.
\end{remark}

In our main results we impose the following assumption on the annealed model.

\begin{assumption}[Metastability of the annealed model] 
  \label{ass:meta}
  For some $(\beta, h)$ $ \in (0,\infty) \times \mathbb{R}$,  the following holds for the Markov chain $(\widetilde{\Sigma}_N(t))_{t \in \mathbb{N}_0}$ of the annealed model. There exist $2 \leq K \in \mathbb{N}$ and $k_1 > 0$ such that, 
  \begin{align}
    \label{eq:def:meta:alt}
    \prob\biggl[
      \liminf_{N \to \infty} \Otimeta(N)
    \biggr]
    =
    1,
  \end{align}
  where $\Otimeta(N)$ is defined in \eqref{eq:def_Otimeta} and depends on $K$ and $k_1$.
\end{assumption}

\begin{remark}
  Assertion of Assumption~\ref{ass:meta} can be rephrased as follows: $\prob$-a.s., there exists a finite random variable $N_0\in \mathbb{N}$ and a sequence $(\{\mathcal{M}_{1,N}, \ldots, \mathcal{M}_{K,N}\})_{N \geq N_0}$ of $K$ non-empty mutually disjoint subsets of $\cS_N$ (possibly depending on $\omega$) such that for any $N \geq N_0$ the process $(\widetilde{\Sigma}_N(t))_{t \in \mathbb{N}_0}$ is $\tilde{\rho}_N \coloneqq \me^{-k_1 N}$-metastable with respect to $\{\mathcal{M}_{1,N}, \ldots, \mathcal{M}_{K,N}\}$.
\end{remark}

\begin{remark}
  Let us illustrate in the case of Example~\ref{ex:dilute:CW} how to identify candidates of metastable sets. It is well know (cf.\ \cite[Section 3.5]{Bo06}) that, for any $\beta \in (0, \infty)$ and $h \in \mathbb{R}$,
  \begin{align}
    \lim_{N \to \infty} \frac{1}{\beta N} \log \widetilde{Z}_N
    =
    -\inf_{x \in [-1,1]} \widetilde{F}_{\beta, h}(x),
  \end{align}
  where $\widetilde{F}_{\beta, h} : [-1,1] \to \mathbb{R}$ denotes the \emph{free energy} per vertex of the annealed model given by
  \begin{align}
    \widetilde{F}_{\beta,h}(x)
    \coloneqq 
    -\frac{1}{2} x^2 - hx
    + \frac{1}{\beta}
    \biggl(
      \frac{1-x}{2} \log \frac{1-x}{2} + \frac{1+x}{2} \log \frac{1+x}{2}
    \biggr).
  \end{align}
  In particular, for any $\beta > \beta_c \coloneqq 1$ and $h \in (-h_c(\beta), h_c(\beta))$, where the critical strength of the magnetic field is given by
  \begin{align*}
    h_c(\beta) 
    \coloneqq 
    \sqrt{1-\beta^2} - \frac{1}{2\beta}
    \log\bigl(\beta (1 + \sqrt{1-1/\beta})^2\bigr),  
  \end{align*}
  the free energy $\widetilde{F}_{\beta, h}$ admits two local minima $m_1, m_2 \in (-1,1)$.  For $N \in \mathbb{N}$, let $m_{1,N}$ and $m_{2,N}$ be the closest point in $\{-1, -1 + 2/N, \ldots, 1 - 2/N, 1\}$ to $m_1$ and $m_2$, respectively. Define the sets $\mathcal{M}_{1,N}, \mathcal{M}_{2,N} \subset \cS_N$ as the (set-valued) pre-image of the empirical magnetization $\cS_N \ni \sigma \mapsto m_N(\sigma) \coloneqq \frac{1}{N} \sum_{i=1}^N \sigma_i$, i.e.,
  \begin{align}
    \mathcal{M}_{1,N} \coloneqq m_N^{-1}(m_{1,N})
    \qquad\text{and}\qquad
    \mathcal{M}_{2,N} \coloneqq m_N^{-1}(m_{2,N}).
  \end{align}
  By using arguments similar to those given in \cite[Lemma~4.1]{SS19}, it follows that $\{\mathcal{M}_{1,N}, \mathcal{M}_{2,N}\}$ forms a pair of metastable sets in the sense of Definition~\ref{def:meta}.
\end{remark}

For fixed $N \in \mathbb{N}$, given the metastable sets $\{\mathcal{M}_{1,N}, \ldots, \mathcal{M}_{K,N}\}$, we can decompose the state space $\cS_N$ into the domains of attraction with to respect the dynamics of the annealed model. More precisely, by following \cite[Definition~1.4]{SS19}, within the event $\Otimeta(N)$ the metastable sets $\{\mathcal{M}_{1,N}, \ldots, \mathcal{M}_{K,N}\}$ give rise to a \emph{metastable partition} $\{\cS_{1,N}, \ldots, \cS_{K,N}\}$ of the state space $\cS_N$ such that
\begin{align}
  \mathcal{M}_{i,N} 
  \subseteq
  \cS_{i,N}
  \subset
  \mathcal{V}_{i,N}, \qquad i \in \{1, \ldots, K\}.
\end{align}
The \emph{local valley} $\mathcal{V}_{i,N}$ around the metastable set $\mathcal{M}_{i,N}$ is defined as
\begin{align}
  \mathcal{V}_{i,N}
  \coloneqq 
  \mathcal{M}_{i,N} 
  \cup
  &\biggl\{
    \sigma \in \mathcal{M}_N^c
    \colon
    \nonumber\\
    &
    \tProb_{\sigma}^{N }\Bigl[
      \tau_{\mathcal{M}_{i,N }}^N
      < \tau_{\mathcal{M}_{N} \setminus \mathcal{M}_{i,N }}^N
    \Bigr]
    \geq
    \max_{j \ne i}
    \tProb_{\sigma}^{N }\Bigl[
      \tau_{\mathcal{M}_{j,N }}^N
      < \tau_{\mathcal{M}_{N} \setminus \mathcal{M}_{j,N }}^N
    \Bigr]
  \biggr\},
\end{align}
where we recall that $\mathcal{M}_N \coloneqq \bigcup_{i=1}^K \mathcal{M}_{i,N}$.

Our first theorem says that, subject to Assumption~\ref{ass:meta}, $(\Sigma_N(t))_{t \in \mathbb{N}_0}$ also exhibits metastable behaviour in the sense of Definition~\ref{def:meta}.

\begin{theorem}[Metastability]
  \label{thm:meta}
  Suppose that $(\beta, h) \in (0, \infty)  \times \mathbb{R}$ satisfies Assumption~\ref{ass:meta}. Then, for any $c_0 \in (0, k_1)$, the event
  \begin{align}
    \Ometa(N)
    \coloneqq 
    \Bigl\{ 
      \omega \in \Otimeta(N)\colon 
      (\Sigma_N(t))_{t \in \mathbb{N}_0}(\omega) \text{ is $\me^{-c_0 N}$-metastable w.r.t.}
      \nonumber\\
      \{\mathcal{M}_{1,N}, \ldots, \mathcal{M}_{K,N}\}(\omega)
    \Bigr\}
    \label{eq:def_Ometa}
  \end{align}
  satisfies 
  \begin{align}\label{eq:meta_thm}
    \prob\biggl[\liminf_{N \to \infty} \Omega_{\textnormal{meta}}(N)\biggr]
    =
    1.
  \end{align}
\end{theorem}

Label the metastable sets $\mathcal{M}_{1,N}, \ldots, \mathcal{M}_{K,N}$ in such a way that they are ordered decreasingly according to their weights under the Gibbs measure of the annealed model, i.e.,  for all $N \in \mathbb{N}$, on the event $\Otimeta(N)$, 
\begin{align}
  \label{eq:ordered:metasets}
  \tilde{\mu}_N\bigl[\mathcal{M}_{1,N}\bigr] 
  \geq
  \tilde{\mu}_N\bigl[\mathcal{M}_{2,N}\bigr]
  \geq
  \ldots
  \geq
  \tilde{\mu}_N\bigl[\mathcal{M}_{K,N}\bigr],
  \qquad \prob\text{-a.s.}
\end{align}

Fix $i \in \{2, \ldots, K\}$ such that there exists a $k_2 \in (0, \infty)$ satisfying 
\begin{align}
  \label{eq:non-degeneracy}
  \tilde{\mu}_N[\mathcal{S}_{j,N}] 
  \leq 
  \me^{-k_2N} \tilde{\mu}_N[\mathcal{S}_{i,N}],
  \qquad \prob\text{-a.s., } \forall j \in \{i+1, \ldots, K\}, \,N \in \mathbb{N},
\end{align}
and set, for any $N \in \mathbb{N}$, $\prob$-a.s.\ on the event $\Otimeta(N)$, 
\begin{align}
  \label{def:A_N:B_N}
  \mathcal{A}_N \coloneqq \mathcal{M}_{i,N},
  \qquad
  \mathcal{B}_N \coloneqq \bigcup_{j=1}^{i-1}\mathcal{M}_{j,N}.
\end{align}
Note that $\cB_N$ is the union of all metastable sets with weight not smaller than the weight of $\cA_N$. Furthermore, note that, since for $i=K$ the condition in \eqref{eq:non-degeneracy}, to which we will also refer as the \emph{non-degeneracy} condition, is void, the existence of such an index $i$ is always guaranteed.  However, if the condition holds for further $i$, one has more freedom in the choice of the set $\cA_N$, making our results more general.

\begin{remark}
  The non-degeneracy condition in \eqref{eq:non-degeneracy} can be relaxed by replacing $\me^{-k_2N}$ with some $ \tilde{\delta}_N$ satisfying $\me^{-k_2N}\leq \tilde{\delta}_N< \me^{-c \sqrt{N}}$ for some sufficiently large $c \in (0, \infty)$. The technical reasons can be found in the proofs in Section~\ref{sec:harm}.
\end{remark}

Before proceeding, we define, for $N \in \mathbb{N}$ and non-empty disjoint sets $\mathcal{A}, \mathcal{B} \subset \cS_N$, the so-called \emph{last-exit biased distribution} on $\mathcal{A}$ for the transition from $\mathcal{A}$ to $\mathcal{B}$ by
\begin{align}
  \label{eq:def_nuAB}
  \nu_{\mathcal{A}, \mathcal{B}}(\sigma)
  \equiv
  \nu_{\cA, \cB}^N(\sigma)
  =
  \frac{\mu_{N}(\sigma) \Prob_{\sigma}^N\bigl[\tau_{\cB}^N < \tau_{\cA}^N\bigr]}
  {
    \sum_{\sigma \in \cA} \mu_{N}(\sigma)
    \Prob_{\sigma}^N\bigl[\tau_{\cB}^N < \tau_{\cA}^N\bigr]
  },
  \qquad \sigma \in \cA.
\end{align}
This distribution plays an essential role in the potential-theoretic approach to meta\-stability, as can be seen in \eqref{eq:PotTheor} below.

We are now ready to state our second theorem, in which we compare the mean hitting time of $\mathcal{B}_N$ for the Markov chain $(\Sigma_N(t))_{t \in \mathbb{N}_0}$ starting from the set $\mathcal{A}_N$ according to the distribution $\nu_{\cA_N, \cB_N}$ with the corresponding quantity for the Markov chain $(\widetilde{\Sigma}_N(t))_{t \in \mathbb{N}_0}$. Under the regular conditional distribution $\prob_{\cG}$, we obtain for the ratio of these metastable hitting times estimates both on its tail behaviour and on its moments. 

\begin{theorem}
  \label{thm_main}
  Suppose that $(\beta, h) \in (0, \infty) \times \mathbb{R}$ satisfies Assumption~\ref{ass:meta}. Set
  \begin{align}
    \label{eq:defAlphaN}
    \alpha_N
    \coloneqq 
    \frac{\beta^2}{2N^2}
    \sum_{\substack{i,j=1 \\ i < j}}^N \Var_{\mathcal{G}}[J_{ij}].
  \end{align}
  \begin{enumerate}[(i)]
  \item 
    For $t \in \mathbb{N}_0$,  $\prob$-a.s.,
    \begin{align}
      \label{eq:main:conc}
      \lim_{N \to \infty}
      \prob_{\cG}\Biggl[
        \me^{-t-\alpha_N}
        \leq
        \frac{
          \Mean_{\nu_{\mathcal{A}_N, \mathcal{B}_N}}^N\bigl[\tau_{\mathcal{B}_N}^N\bigr]
        }{
          \tMean_{\tilde{\nu}_{\mathcal{A}_N, \mathcal{B}_N}}^N\bigl[\tilde{\tau}_{\mathcal{B}_N}^N\bigr]
        }
        \leq
        \me^{+t+2\alpha_N}
      \Biggr]
      \geq
      1 - 4\, \me^{-t^2 / (2 \beta k_J)^2}.
    \end{align}
  \item For any $q \geq 1$ and $c\in (0, \infty)$, let
    \begin{align}
      &\Omega_{q,c}(N) 
      \coloneqq 
      \nonumber \\
      &
      \Biggl\{ 
        \omega \colon \me^{-\alpha_N}\Bigl(1 - \tfrac{c}{\sqrt{N}}\Bigr)
        \leq
        \frac{
          \mean_{\cG}\Bigl[
            \Mean_{\nu_{\mathcal{A}_N, \mathcal{B}_N}}^N\bigl[\tau_{\mathcal{B}_N}^N\bigr]^q
          \Bigr]^{1/q}(\omega)
        }
        {
          \tMean_{\tilde{\nu}_{\mathcal{A}_N, \mathcal{B}_N}}^N\bigl[\tilde{\tau}_{\mathcal{B}_N}^N\bigr](\omega)
        }
        \leq
        \me^{4q \alpha_N}\Bigl(1 + \tfrac{c}{\sqrt{N}}\Bigr)
      \Biggr\}.
      \label{eq:def_Oqc}
    \end{align}
    Then, for any $q \geq 1$ there exists $c_1 \in (0, \infty)$ such that 
    \begin{align}
      \label{eq:main:moments}
      \prob\Biggl[\liminf_{N \to \infty} \Omega_{q,c_1}(N)\Biggr] = 1 .
    \end{align}
  \end{enumerate}
\end{theorem}

\begin{remark}
  \begin{enumerate}[(a)]
  \item 
    Since the random variables $(J_{ij})_{1 \leq i < j < \infty}$ are assumed to be uniformly bounded, it follows that $\alpha_N=O(1)$.
    
  \item 
    The non-degeneracy condition \eqref{eq:non-degeneracy} embedded in the definition \eqref{def:A_N:B_N} of $\cA_N$ ensures that, if $i<K$, the metastable sets $\cM_{j,N}$, $j \in \{i+1,\dots, K\}$, are not relevant for the analysis of the crossing times from $\cA_N$ to $\cB_N$.

  \item      The choice of the initial configuration being drawn from the (quenched) last-exit biased distribution plays an important role in the potential-theoretic approach to metastability.  Indeed, under this particular initial distribution supported on (the inner boundary of) the set $\mathcal{A}_{N}$ the mean hitting time, $\Mean_{\sigma}\bigl[\tau_{\mathcal{B}_{N}}^{N}\bigr]$, of the set $\mathcal{B}_{N}$ can be represented in terms of two analytic quantities: the harmonic function, $h_{\mathcal{A}_{N}, \mathcal{B}_{N}}$, and the capacity, $\capa_N(\mathcal{A}_{N}, \mathcal{B}_{N})$, see \eqref{eq:PotTheor}.  If $\mathcal{A}_{N}$ and  $\mathcal{B}_{N}$ are metastable sets, one expects that the function $\sigma \mapsto \Mean_{\sigma}\bigl[\tau_{\mathcal{B}_{N}}^{N}\bigr]$ is almost constant on the set $\mathcal{A}_{N}$.  Such regularity statement has been proven for certain spin systems by means of coupling methods, cf.\ \cite{BBI12, MOS90, MS88}.  The construction of a suitable coupling in the present setting, while expected to be possible, requires some highly non-trivial adaptation of the previously used strategies, due to the random interaction, that goes beyond the scope of the present paper.

  \end{enumerate}
\end{remark}

\begin{remark}
  \begin{enumerate}[(a)]
  \item      The results still hold true for different reversible, single spin update dynamics.  However, the Glauber dynamics with Metropolis transition rate is a standard choice.

  \item      The assertions of both Theorem~\ref{thm:meta} and \ref{thm_main} are still valid beyond the regime of uniformly bounded random variables.  The forthcoming paper \cite{DLPS24+} studies the Potts model with random interactions assuming that $(J_{ij})_{1 \leq i < j < \infty}$ are i.i.d.\ unbounded random variables with finite exponential moments.

  \end{enumerate}
\end{remark}


\subsection{Methods and outline}
\label{sec:meth}


\subsubsection{Key notions from the potential-theoretic approach to metastability}

To prove Theorems~\ref{thm:meta} and \ref{thm_main}, we crucially rely on  potential theory, which allows us to express probabilistic objects of interest in terms of solutions of certain boundary value problems. It is well-known (see e.g.\ \cite[Corollary~7.11]{BdH15}) that, for any $N \in \mathbb{N}$ and any non-empty disjoint $\mathcal{A}, \mathcal{B} \subset \cS_N$, the mean hitting time of $\mathcal{B}$ starting from the last-exit biased distribution $\nu_{\mathcal{A}, \mathcal{B}}^N$ on $\mathcal{A}$ is given by
\begin{align}
  \label{eq:PotTheor}
  \Mean_{\nu_{\mathcal{A}, \mathcal{B}}}^N\bigl[\tau_{\mathcal{B}}^N\bigr]
  = \frac{\norm{h_{\mathcal{A}, \mathcal{B}}^N}_{\mu_N}}
  {\capa_N(\mathcal{A}, \mathcal{B})},
\end{align}
where $\norm{h_{\mathcal{A}, \mathcal{B}}^N}_{\mu_N}$ denotes the $\ell^1(\mu_N)$-norm of the  \emph{equilibrium potential} $h_{\mathcal{A}, \mathcal{B}}^N$ of the pair $(\mathcal{A}, \mathcal{B})$, i.e., the function $h_{\mathcal{A}, \mathcal{B}}^N : \cS_N \to [0,1]$ that is the unique solution of the boundary value problem
\begin{align}
  \left\{
    \begin{array}{rcll}
      \big(L_N f\big)(\sigma)
      &\mspace{-5mu}=\mspace{-5mu}& 0,
      &\sigma \in \cS_N \setminus (\mathcal{A} \cup \mathcal{B}),\\ 
      f(\sigma)
      &\mspace{-5mu}=\mspace{-5mu}& \ind_{\mathcal{A}}(\sigma), \quad 
      &\sigma \in \mathcal{A} \cup \mathcal{B}.
    \end{array}
  \right.
\end{align}
Note that the equilibrium potential has a natural interpretation in terms of \emph{hitting probabilities}, namely, $h_{\mathcal{A}, \mathcal{B}}^N(\sigma) = \Prob_{\sigma}^N\bigl[\tau_{\mathcal{A}}^N < \tau_{\mathcal{B}}^N\bigr]$, for all $\sigma \in \cS_N \setminus (\mathcal{A} \cup \mathcal{B})$. The \emph{capacity} $\capa_N(\mathcal{A}, \mathcal{B})$ of the pair $(\mathcal{A}, \mathcal{B})$ is defined by
\begin{align}
  \label{eq:def_capa}
  \capa_N(\mathcal{A}, \mathcal{B})
  \coloneqq 
  \sum_{\sigma \in \mathcal{A}} \mu_N(\sigma) 
  \Prob_{\sigma}^N\bigl[\tau_{\mathcal{B}}^N < \tau_{\mathcal{A}}^N\bigr]
  =
  \sum_{\sigma \in \mathcal{A}} \mu_N(\sigma)
  \bigl(-L_N h_{\mathcal{A}, \mathcal{B}}^N\bigr)(\sigma).
\end{align}
From this definition it is clear that
\begin{align}
  \label{eq:P_hit_cond}
  \Prob_{\mu_N \vert \mathcal{A}}^N\bigl[\tau_{\mathcal{B}}^N < \tau_{\mathcal{A}}^N \bigr] 
  = \frac{\capa_N(\mathcal{A}, \mathcal{B})}{\mu_N[\cA]},
\end{align}
where $\mu_N \vert \cA$ denotes the Gibbs measure $\mu_N$ conditioned on the set $\cA$.

Capacities can be expressed in terms of variational principles that are very useful to obtain upper und lower bounds (see \cite[Section~7.3]{BdH15} for more details).  Upper bounds are obtained by using the \emph{Dirichlet principle}, which states that
\begin{align}
  \label{eq:Dirichlet}
  \capa_N(\mathcal{A}, \mathcal{B})
  =
  \inf \bigl\{
    \mathcal{E}_N(f) \colon f \in \mathcal{H}_{\mathcal{A}, \mathcal{B}}^N
  \bigr\}.
\end{align}
Here, $\mathcal{H}_{\mathcal{A}, \mathcal{B}}^N$ denotes the set of all functions from $\cS_N$ to $\mathbb{R}$ that are equal to $1$ on $\mathcal{A}$ and $0$ on $\mathcal{B}$, and
\begin{align}
  \mathcal{E}_N(f)
  \coloneqq 
  \scpr{f, -L_N f}_{\mu_N}
  =
  \frac{1}{2} \sum_{\sigma, \sigma' \in \cS_N}
  \mu_N(\sigma) p_N(\sigma, \sigma') \bigl(f(\sigma) - f(\sigma')\bigr)^2
\end{align}
is the \emph{Dirichlet form}. We recall that the transition probabilities $p_N$ are defined in \eqref{eq:defpN}.

Lower bounds are obtained via the \emph{Thomson principle}, which states that
\begin{align}
  \label{eq:Thomson}
  \capa_N(\mathcal{A}, \mathcal{B})
  =
  \sup\Bigl\{
    \frac{1}{\mathcal{D}_N(\varphi)} \colon
    \varphi \in \mathcal{U}_{\mathcal{A}, \mathcal{B}}^N
  \Bigr\}
  =
  \Bigl(
    \inf\bigl\{
      \mathcal{D}_N(\varphi)
      \colon
      \varphi \in \mathcal{U}_{\mathcal{A}, \mathcal{B}}^N
    \bigr\}
  \Bigr)^{-1},
\end{align}
where $\mathcal U_{\mathcal{A},\mathcal{B}}^N$ is the space of all unit antisymmetric $\mathcal{A}\mathcal{B}$-flows $\varphi: \cS_N \times \cS_N \to \mathbb{R}$, and
\begin{align}
  \mathcal{D}_N(\varphi)
  \coloneqq 
  \frac{1}{2} 
  \sum_{\substack{\sigma, \sigma' \in \cS_N \\ \sigma \sim \sigma'}}
  \frac{\varphi(\sigma, \sigma')^2}{\mu_N(\sigma)\, p_N(\sigma, \sigma')}.
\end{align}

\subsubsection{Strategy of proofs}

The proof of Theorem~\ref{thm:meta} relies on Definition~\ref{def:meta} and on \eqref{eq:P_hit_cond}, together with an application of the Dirichlet principle in combination with a comparison of the quenched Hamiltonian $H_N$ and the annealed Hamiltonian $\widetilde{H}_N$ on a particular event of high probability.

We prove Theorem~\ref{thm_main}(i) by combining concentration inequalities for the logarithm of the mean hitting time $\Mean_{\nu_{\mathcal{A}_N, \mathcal{B}_N}}^N\bigl[\tau_{\mathcal{B}_N}^N\bigr]$ of the quenched model with bounds on the distance between the (conditional on $\cG$) mean of that logarithm and the logarithm of the mean hitting time $\tMean_{\tilde{\nu}_{\mathcal{A}_N, \mathcal{B}_N}}^N\bigl[\tilde{\tau}_{\mathcal{B}_N}^N\bigr]$ of the annealed model. Estimates of the latter type, comparing conditional means with average means, will be called \emph{annealed estimates}.  The results in Theorem~\ref{thm_main}(ii) are annealed estimates as well. In view of \eqref{eq:PotTheor}, estimates on the mean hitting time $\Mean_{\nu_{\mathcal{A}_N, \mathcal{B}_N}}^N[\tau_{\mathcal{B}_N}^N]$, or on its logarithm, will follow once we have separately proven corresponding estimates for both $Z_N \capa_N(\mathcal{A}_N, \mathcal{B}_N)$ and $Z_N \norm{h_{\mathcal{A}_N, \mathcal{B}_N}^N}_{\mu_N}$. 

To prove concentration inequalities for the quantities $\log[ Z_N \capa_N(\mathcal{A}_N, \mathcal{B}_N)]$ and $\log [Z_N \norm{h_{\mathcal{A}_N, \mathcal{B}_N}^N}_{\mu_N}]$, we use a conditional version of McDiarmid's bounded differences inequality (see Proposition~\ref{prop:McDiarmid} below). This strategy for proving concentration is different from the one used in \cite{BMP21}, where Talagrand's concentration inequality was used. The advantage of McDiarmid over Talagrand is twofold. First, McDiarmid's inequality provides exact constants. Second, it does not require convexity of the map $J \mapsto \log \bigl(Z_N \capa_N(\mathcal{A}_N, \mathcal{B}_N)\bigr)$, which is crucial because we do not know how to prove convexity.  

Estimates  on capacities for Theorem~\ref{thm_main} are proven by using the Dirichlet principle and the Thomson principle, and do not require any assumption on metastability. Finding estimates on the equilibrium potential, however, is more involved. We use  a result that is similar to \cite[Theorem~1.7]{SS19} (Proposition~\ref{prop:harm:estimate} below), for which the non-degeneracy assumption in \eqref{eq:non-degeneracy} is required, together with the same comparison of the Hamiltonians $H_N$ and $\widetilde{H}_N$ that is used in the proof of Theorem~\ref{thm:meta}, both holding with high probability. We emphasise that the constants appearing in our statements may depend on the parameters of the model.


\subsubsection{Outline}

The remainder of the paper is organised as follows. In Section~\ref{sec:meta_proof} we provide the proof of Theorem~\ref{thm:meta} on metastability of the quenched model. In Section~\ref{sec:capa} we provide estimates on capacities. Section~\ref{sec:harm} is devoted to stating and proving estimates on weighted sums of the equilibrium potential, called harmonic sums. In Section~\ref{sec:mean:hitting:times} we prove Theorem~\ref{thm_main} by using the results of the previous sections. Appendix~\ref{app:McD} contains the conditional version of the McDiarmid's inequality that is used in the paper.

\section{Metastability}
\label{sec:meta_proof}

Before proving Theorem~\ref{thm:meta}, we provide in Lemma~\ref{lemma:conc:H_N} a comparison of the quenched Hamiltonian $H_N$ and the annealed Hamiltonian $\widetilde{H}_N$. This lemma will be used both in the proof of Theorem~\ref{thm:meta} below and in Section~\ref{sec:harm}, where we deal with estimates on the equilibrium potential.

Given a positive real sequence $(a_N)_{N \in \mathbb{N}}$, let
\begin{align}
  \label{eq:def:Xi(N)}
  \Xi(a_N)
  \coloneqq 
  \Bigl\{
    \max_{\sigma \in \cS_N} \bigl\vert H_N(\sigma) - \widetilde{H}_N(\sigma) \bigr\vert
    <
    a_N
  \Bigr\}
  \subset \Omega, \qquad N \in \mathbb{N} \setminus\{1\},
\end{align}
denote the event that, uniformly in $\sigma \in \cS_N$, $H_N$ differs from $\widetilde{H}_N$ by at most $a_N$. On the event $\Xi(a_N)$ we have control on the difference between the quantities determining $(\Sigma_N(t))_{t \in \mathbb{N}_0}$ and $(\widetilde{\Sigma}_N(t))_{t \in \mathbb{N}_0}$. Moreover, for suitably chosen sequences $(a_N)_{N \in \mathbb{N}}$, the event $\Xi(a_N)^c$ turns out to be negligible in the limit as $N\to\infty$. 

\begin{lemma}
  \label{lemma:conc:H_N}
  For a positive real sequence $(a_N)_{N \in \mathbb{N}}$, set $b_N \coloneqq a_N^2 / k_J^2 - (N+1) \log 2$. Then, $\prob$-a.s.,
  \begin{align}
    \label{eq:conc:H_N}
    \prob_{\mathcal{G}}\bigl[\Xi(a_N)^c
    \bigr]
    \leq
    \me^{-b_N} \wedge 1,
    \qquad\forall N \in \mathbb{N} \setminus \{1\}.
  \end{align}
\end{lemma} 

\begin{proof}
  Fix $2 \leq N \in \mathbb{N}$. Clearly, it suffices to prove \eqref{eq:conc:H_N} in case $b_N > 0$. 
  Notice that the map $J \longmapsto H_N(\sigma)$ satisfies a bounded difference as \eqref{eq:bounded:difference} with constant $2 k_J/N$ uniformly in $\sigma \in \cS_N$, because $ H_N(\sigma)$ depends linearly on the random coupling $(J_{ij})_{1 \leq i < j \leq N}$ and we assumed in Section~\ref{sec:model} that the $\vert J_{ij}\vert  \leq k_J$ $\prob$-a.s.\ for all $1 \leq i < j < \infty$. 
  Since in addition the triangular array $J=(J_{ij})_{1 \leq i < j < \infty}$ is assumed to be conditionally independent given $\mathcal{G}$, we can apply McDiarmid's concentration inequality (Proposition~\ref{prop:McDiarmid}) with $v=k_J^2(N-1)/(2N)$, together with a union bound, to get that, $\prob$-a.s.,
  \begin{align}
    \begin{split}
      \prob_{\mathcal{G}}\Bigl[
        \max_{\sigma \in \cS_N} \bigl\vert H_N(\sigma) - \widetilde{H}_N(\sigma)\bigr\vert
        \geq
        a_N
      \Bigr]
      &\leq
      \sum_{\sigma \in \cS_N}
      \prob_{\mathcal{G}} \Bigl[
        \bigl\vert H_N(\sigma) - \widetilde{H}_N(\sigma)\bigr\vert \geq a_N
      \Bigr]
      \\
      &\leq
      2^{N+1} \exp\biggl(-\frac{a_N^2 N}{k_J^2 (N-1)}\biggr),
    \end{split}
  \end{align}
  where the additional factor $2$ comes from the absolute value. 
  Since $ (N-1)/N \leq 1$, \eqref{eq:conc:H_N} follows.
\end{proof}

\begin{proof}[Proof of Theorem~{\upshape\ref{thm:meta}}]
  We will prove that
  \begin{align}
    \label{eq:supThm}
    \prob\biggl[\limsup_{N \to \infty} \Omega_{\textnormal{meta}}(N)^c\biggr]
    =
    0,
  \end{align}
  which is equivalent to \eqref{eq:meta_thm}.
  First note that, by the choice of $\mu_N$ and $p_N$ in \eqref{eq:def_muN} and \eqref{eq:defpN},  
  \begin{align}
    Z_N \mu_N(\sigma) \, p_N(\sigma, \sigma')
    = \frac{1}{N}
    \me^{-\beta (H_N(\sigma) \vee H_N(\sigma'))},
    \qquad \sigma \sim \sigma' \in \cS_{N}.
  \end{align}
  An elementary computation yields that, on the event $\Xi(a_N)$, 
  \begin{align}
    \widetilde{H}_N(\sigma) \vee \widetilde{H}_N(\sigma') - a_N
    \leq
    H_N(\sigma) \vee H_N(\sigma')
    \leq
    \widetilde{H}_N(\sigma) \vee \widetilde{H}_N(\sigma') + a_N, \quad \sigma,\sigma' \in \cS_N. 
  \end{align}
  Thus, by a comparison of Dirichlet forms it follows that
  \begin{align}
    \widetilde{Z}_N \widetilde{\mathcal{E}}_N(f)\, \me^{-\beta a_N}
    \leq
    Z_N \mathcal{E}_N(f)
    \leq
    \widetilde{Z}_N \widetilde{\mathcal{E}}_N(f)\, \me^{\beta a_N},
  \end{align}
  for any $f : \cS_N \to \mathbb{R}$. In view of the Dirichlet principle \eqref{eq:Dirichlet}, we deduce that, on the event $\Xi(a_N)$, 
  \begin{align}
    \label{eq:comparison:cap}
    \me^{-\beta a_N}
    \leq
    \frac{Z_N \capa_N(\mathcal{X}, \mathcal{Y})}
    {\widetilde{Z}_N \widetilde{\capa}_N(\mathcal{X}, \mathcal{Y})}
    \leq
    \me^{\beta a_N},
    \qquad \emptyset \neq \mathcal{X}, \mathcal{Y} \subset \cS_N \text{ disjoint}.
  \end{align}
  Moreover, for any $2 \leq N \in \mathbb{N}$, on the event $\Xi(a_N)$,
  \begin{align}
    \label{eq:comparison:mu}
    \me^{-\beta a_N}
    \leq
    \frac{Z_N \mu_N\bigl[\mathcal{X}\bigr]}
    {\widetilde{Z}_N \tilde{\mu}_N\bigl[\mathcal{X}\bigr]}
    \leq
    \me^{\beta a_N},
    \qquad \emptyset \neq \mathcal{X} \subset \cS_N.
  \end{align}
  It follows from \eqref{eq:P_hit_cond}, \eqref{eq:comparison:cap} and \eqref{eq:comparison:mu} that, on the event $\Xi(a_N)$,
  \begin{align}
    \label{eq:comparison:exit_prob}
    \me^{-2\beta a_N}
    \leq
    \frac{
      \Prob_{\mu_N \vert \mathcal{X}}^N\bigl[
        \tau_{\mathcal{Y}}^N < \tau_{\mathcal{X}}^N
      \bigr]
    }{
      \tProb_{\mu_N \vert \mathcal{X}}^N\bigl[
        \tilde{\tau}_{\mathcal{Y}}^N < \tilde{\tau}_{\mathcal{X}}^N
      \bigr]
    }
    \leq
    \me^{2\beta a_N},
    \qquad \emptyset \neq\mathcal{X}, \mathcal{Y} \subset \cS_N \text{ disjoint}.
  \end{align}
  Thus, on the event $\Xi(a_N) \cap \Otimeta(N)$,
  \begin{align}
    \max_{j \in \{1, \ldots, K\}}&
    \Prob_{\mu_N \vert \mathcal{M}_{j,N}}^N\Bigl[
      \tau_{\mathcal{M}_N \setminus \mathcal{M}_{j,N}}^N < \tau_{ \mathcal{M}_{j,N}}^N
    \Bigr]
    \nonumber \\
    &\leq
    \frac{\tilde{\rho}_N}{K}\, \me^{4 \beta a_N}
    \min_{\mathcal{X} \subset \cS_N \setminus \mathcal{M}_N} 
    \Prob_{\mu_N \vert \mathcal{X}}^N\Bigl[
      \tau_{\mathcal{M}_N}^N < \tau_{\mathcal{X}}^N
    \Bigr].
    \label{eq:meta_ann}
  \end{align}
  Now set $a_N = k_J\sqrt{N k_1+ (N+1)\log 2}$ for $2 \leq N \in \mathbb{N}$, and note that, with this choice of $a_N$, Lemma~\ref{lemma:conc:H_N} implies that
  \begin{align}
    \label{eq:tail:estimate}
    \prob\bigl[\Xi(a_N)^c\bigr] 
    = 
    \mean\bigl[\prob_{\mathcal{G}}\bigl[\Xi(a_N)^c\bigr]\bigr]
    \leq
    \me^{-b_N}
    =
    \me^{-k_1 N},
    \qquad N \in \mathbb{N} \setminus\{1\}. 
  \end{align}
  Recall that $c_0$ is the constant in the definition of $ \Ometa(N)$ in \eqref{eq:def_Ometa}, in Theorem~\ref{thm:meta}. By choosing $N(k_1, c_0, \beta, k_J) \in \mathbb{N}$ in such a way that $c_0 < k_1-4\beta a_N/N$ for all $N \geq N(k_1, c_0, \beta, k_J)$, it follows from \eqref{eq:meta_ann} that $\Xi(a_N) \cap \Otimeta(N) \subseteq \Ometa(N)$ for all $N \geq N(k_1, c_0, \beta, k_J)$. In particular,
  \begin{align}
    \label{eq:OmetaC}
    \Ometa(N)^c 
    \subseteq 
    \Xi(a_N)^c \cup \Otimeta(N)^c, 
    \qquad N \geq N(k_1, c_0, \beta, k_J).
  \end{align}
  Therefore, using continuity of the probability measure, we get
  \begin{align}
    \prob\biggl[\limsup_{N \to \infty} \Ometa(N)^c\biggr]\;
    &\overset{\mspace{-15mu}\eqref{eq:OmetaC}\mspace{-15mu}}{\leq}
    \;\prob\biggl[
      \limsup_{N \to \infty}
      \Bigl( \Xi(a_N)^c \cup \Otimeta(N)^c \Bigr)
    \biggr]
    \nonumber\\
    &\leq
    \lim_{N \to \infty} 
    \biggl(
      \prob\biggl[\bigcup\nolimits_{m \geq N} \Xi(a_m)^c \biggr] 
      +
      \prob\biggl[
        \bigcup\nolimits_{m \geq N} \Otimeta(m)^{c}
      \biggr]
    \biggr)
    \nonumber\\
    &=
    \prob\biggl[ \limsup_{N \to \infty} \Xi(a_N)^c \biggr] 
    +
    \prob\biggl[
      \limsup_{N \to \infty} \Otimeta(N)^{c}
    \biggr].
  \end{align}
  Since, by \eqref{eq:tail:estimate}, $\sum_{N=1}^{\infty} \prob\bigl[\Xi(a_N)^c\bigr] < \infty$, an application of the Borel Cantelli Lemma and Assumption~\ref{ass:meta} yields that
  \begin{align}
    \prob\biggl[\limsup_{N \to \infty} \Ometa(N)^c\biggr]
    \leq
    \prob\biggl[ \limsup_{N \to \infty} \Xi(a_N)^c \biggr] 
    +
    \prob\biggl[
      \limsup_{N \to \infty} \Otimeta(N)^{c}
    \biggr]
    =
    0.
  \end{align}
\end{proof}


\section{Capacity estimates}
\label{sec:capa}

In this section we provide general estimates on the capacity of the quenched model compared to the annealed model. These estimates are general because we do not require any assumption on metastability, and the sets involved in the estimates are general disjoint subsets of the configuration space.

In Section~\ref{sec:cap_conc} we prove concentration for the logarithm of the capacities by using the Dirichlet principle and McDiarmid's concentration inequality. In Section~\ref{sec:annealed:cap} we first estimate the conditional mean of $Z_N \mu_N$ and $p_N$ in terms of the corresponding quantities of the annealed model $\widetilde{Z}_N \tilde{\mu}_N$ and $\tilde{p}_N$ in Lemma~\ref{lem:EB:Delta+HvH}, and afterwards prove annealed capacity estimates by using both the Dirichlet and the Thomson principle, together with Lemma~\ref{lem:EB:Delta+HvH}. The latter is crucial also in the proof of annealed estimates of $\norm{h_{\mathcal{A}_N,\mathcal{B}_N}^N}_{\mu_N}$ in Section~\ref{sec:h_ann}.

All formulas in this Section~\ref{sec:capa} are intended to hold $\prob$-a.s. In order to lighten notation, we refrain from repeating that.


\subsection{Concentration of quenched capacities} \label{sec:cap_conc}

\begin{proposition}
  \label{pro:conc_logcap}
  Let $2 \leq N \in \mathbb{N}$, and consider two non-empty disjoint subsets $\mathcal{X}, \mathcal{Y}  \subset \cS_N$. Then, for any $t \in \mathbb{N}_0$,
  \begin{align}
    \label{eq:conc_logcap}
    \prob_{\mathcal{G}}\Bigl[
      \bigl\vert
        \log \bigl(Z_N \capa_N(\mathcal{X}, \mathcal{Y})\bigr)
        -
        \mean_{\mathcal{G}}\bigl[
          \log\bigl(Z_N \capa_N(\mathcal{X}, \mathcal{Y})\bigr)
        \bigr]
      \bigr\vert
      > t
    \Bigr]
    \leq
    2\, \me^{-t^2/(\beta k_J)^2}.
  \end{align}
\end{proposition}

\begin{proof}
  First, recall that the triangular array $(J_{ij})_{1 \leq i < j < \infty}$ is assumed to be conditionally independent given $\mathcal{G}$. Hence, in view of McDiarmid's concentration inequality (Proposition~\ref{prop:McDiarmid}), the assertion in \eqref{eq:conc_logcap} is immediate once we show that, for any $2 \leq N \in \mathbb{N}$, the mapping
  \begin{align}
    (J_{ij})_{1 \leq i < j \leq N} 
    \longmapsto 
    F_N\bigl((J_{ij})_{1 \leq i < j \leq N}\bigr) 
    \coloneqq 
    \log \bigl(Z_N \capa_N(\mathcal{X}, \mathcal{Y})\bigr)
  \end{align}
  satisfies a bounded difference estimate. More precisely, it is sufficient to show that, for any $1 \leq k < l \leq N$,
  \begin{align}
    \label{eq:bounded:diff:cap}
    \bigl\vert
      F_N\bigl((J_{ij})_{1 \leq i < j \leq N}\bigr)
      - F_N\bigl((J_{ij}')_{1 \leq i < j \leq N}\bigr)
    \bigr\vert
    \leq
    \frac{2 \beta k_J}{N},
  \end{align}
  where $J_{ij}' \coloneqq J_{ij}$ for all $1 \leq i < j \leq N$ such that $(i,j) \ne (k,l)$, and $J_{kl}'$ is a conditionally independent copy of $J_{kl}$ given $\mathcal{G}$. In the sequel, we write $H_N^J$, $Z_N^J$, $\mathcal{E}_N^{J}$ and $\capa_N^J(\mathcal{X}, \mathcal{Y})$ to emphasise the dependence on the random coupling $J=(J_{ij})_{1 \leq i < j \leq N}$.  
  
  We proceed by following the same line of argument that led to \eqref{eq:comparison:cap} in the proof of Theorem~\ref{thm:meta}. By 
  \begin{align}
    \label{eq:H_N:JvsJ'} 
    \bigl\vert H_N^J(\sigma) - H_N^{J'}(\sigma) \bigr\vert 
    = 
    \frac{\bigl\vert J_{kl} - J_{kl}'\bigr\vert}{N} 
    \leq 
    \frac{2k_J}{N},
    \qquad \sigma \in \cS_N,
  \end{align}
  an elementary computation yields that, for any $\sigma, \sigma' \in \cS_N$,
  \begin{align}
    H_N^{J'}(\sigma) \vee H_N^{J'}(\sigma') - \frac{2k_J}{N}
    \leq
    H_N^{J}(\sigma) \vee H_N^{J}(\sigma')
    \leq
    H_N^{J'}(\sigma) \vee H_N^{J'}(\sigma') + \frac{2k_J}{N}.
  \end{align}
  Thus, by a comparison of Dirichlet forms, we obtain, for any $f : \cS_N \to \mathbb{R}$,
  \begin{align}
    Z_N^{J'} \mathcal{E}_N^{J'}(f)\, \me^{-2\beta k_J /N}
    \leq
    Z_N^{J} \mathcal{E}_N^{J}(f)
    \leq
    Z_N^{J'} \mathcal{E}_N^{J'}(f)\, \me^{2\beta k_J /N}.
  \end{align}
  In view of the Dirichlet principle \eqref{eq:Dirichlet}, we deduce that
  \begin{align}
    Z_N^{J'} \capa_N^{J'}(\mathcal{X}, \mathcal{Y})\, \me^{-2 \beta k_J / N}
    \leq
    Z_N^J \capa_N^J(\mathcal{X}, \mathcal{Y})
    \leq
    Z_N^{J'} \capa_N^{J'}(\mathcal{X}, \mathcal{Y})\, \me^{2 \beta k_J / N},
  \end{align}
  which yields \eqref{eq:bounded:diff:cap}.
\end{proof}


\subsection{Annealed capacity estimates} \label{sec:annealed:cap}

\begin{notation}
  For any three sequences $(a_N)_{N \geq 0}, (b_N)_{N \geq 0}$, $ (c_N)_{N \geq 0}$ and $N \in \mathbb{N}$, the notation $a_N = b_N +O(c_N)$ means that there exists a $C \in(0, \infty)$ independent of $\omega$ and $N$ such that 
  \begin{align}
    - C  c_N \leq a_N - b_N \leq C  c_N.
  \end{align}
\end{notation}

Before proving annealed capacity estimates, we prove the following lemma which is used both in the current section and for proving further annealed estimates in Section~\ref{sec:h_ann}.
\begin{lemma}
  \label{lem:EB:Delta+HvH}
  For $2 \leq N \in \mathbb{N}$ the following hold:
  \begin{enumerate}[(i)]
  \item For any $\sigma \in \cS_N$, 
    \begin{align}
      \label{eq:EB:Delta}
      \mean_{\mathcal{G}} \Bigl[\me^{\pm \beta \Delta_N(\sigma)}\Bigr]
      =
      \me^{\alpha_N} \bigl(1 + O(N^{-1})\bigr),
    \end{align}
  \item For any $\sigma, \sigma' \in \cS_N$ with $\sigma \sim \sigma'$,
    \begin{align}
      \label{eq:EB:HvH}
      \mean_\mathcal{G}\Bigl[\me^{\pm \beta (H_N(\sigma) \vee H_N(\sigma'))}\Bigr]
      =
      \me^{\pm \beta(\widetilde{H}_N(\sigma) \vee \widetilde{H}_N(\sigma'))}\, 
      \me^{\alpha_N}
      \bigl(1 + O(N^{-1/2})\bigr),
    \end{align}
  \end{enumerate}
  where $\alpha_N$ is defined in \eqref{eq:defAlphaN} and $\Delta_N(\sigma) \coloneqq H_N(\sigma) - \widetilde{H}_N(\sigma)$, $\sigma \in \cS_N$.
\end{lemma}

\begin{proof}
  \textit{(i)} Denote by $\mathbb{R} \ni t \mapsto \Lambda_{ij}(t) \coloneqq \log \mean_{\mathcal{G}}\bigl[\exp\bigl(t(J_{ij} - \mean_{\mathcal{G}}[J_{ij}])\bigr)\bigr]$ the conditional log-moment generating function given $\mathcal{G}$. By a Taylor expansion up to the third order, we get that, for any $t \in \mathbb{R}$,
  \begin{align}
    \Lambda_{ij}(t)
    =
    \frac{t^2}{2} \Var_{\mathcal{G}}\bigl[J_{ij}\bigr]
    + \frac{t^3}{2}
    \int_0^1 (1-\theta)^2 \Lambda_{ij}'''(\theta t) \md \theta.
  \end{align}
  Since the random variables are assumed to be uniformly bounded, i.e., $\vert J_{ij}\vert  \leq k_J$, an elementary computation exploiting Cram\'{e}r's measure yields that, $\vert \Lambda_{ij}'''(t)\vert  \leq 6k_J^3$. Hence
  \begin{align}
    \bigl\vert \Lambda_{ij}(t) - \frac{t^2}{2} \Var_{\mathcal{G}}[J_{ij}] \bigr\vert
    \leq
    k_J^3 \abs{t}^3.
  \end{align}
  Since the triangular array $(J_{ij})_{1 \leq i < j < \infty}$ is conditionally independent given $\mathcal{G}$, we have
  \begin{align}
    \Bigl\vert
      \log \mean_{\mathcal{G}}\Bigl[\me^{\pm \beta \Delta_N(\sigma)}\Bigr]
      -
      \alpha_N
    \Bigr\vert
    \leq
    \sum_{\substack{i,j=1 \\ i < j}}^N
    \Bigl\vert
      \Lambda_{ij}\bigl(\tfrac{\mp\beta}{N} \sigma_i \sigma_j\bigr)
      - \frac{(\mp\beta \sigma_i \sigma_j)^2}{2N^2} 
      \Var_{\mathcal{G}}\bigl[J_{ij}\bigr]
    \Bigr\vert
    \leq
    \frac{(\beta k_J)^3}{2N},
  \end{align}
  which concludes the proof of \eqref{eq:EB:Delta}. In particular, for any $\sigma \in \cS_N$,
  \begin{align}
    \label{eq:EB:H}
    \me^{\pm \beta \widetilde{H}_N(\sigma) + \alpha_N}\, \me^{-(\beta k_J)^3/2N}
    \leq
    \mean_{\mathcal{G}}\Bigl[\me^{\pm \beta H_N(\sigma)}\Bigr]
    \leq
    \me^{\pm \beta \widetilde{H}_N(\sigma) + \alpha_N}\, \me^{(\beta k_J)^3/2N}. 
  \end{align}
  
  \noindent  
  \textit{(ii)} Because the proofs of \eqref{eq:EB:HvH} for $\pm \beta$ are similar, we give a detailed proof for $+\beta$ only.  Since the conditional expectation of the maximum of two random variables is bounded from below by the maximum of their conditional expectations, it is immediate from \eqref{eq:EB:H} that, for any $\sigma, \sigma' \in \cS_N$,
  \begin{align}
    \mean_{\mathcal{G}}\Bigl[\me^{\beta(H_N(\sigma) \vee H_N(\sigma'))}\Bigr]
    \geq
    \me^{\beta (\widetilde{H}_N(\sigma) \vee \widetilde{H}_N(\sigma'))}\,
    \me^{\alpha_N} \Bigl(1 - \frac{(\beta k_J)^3}{2N}\Bigr).
  \end{align}
  Thus, we are left with proving the desired upper bound. For this purpose, we define, for $\sigma \in \cS_N$ and $k \in \{1, \ldots, N\}$,
  \begin{align}
    H_N^k(\sigma)
    \coloneqq 
    -\sigma_k
    \biggl(
      \frac{1}{N} \sum_{\substack{j=1\\j > k}}^N J_{kj} \sigma_j
      + \frac{1}{N} \sum_{\substack{i=1\\i<k}}^N J_{ik} \sigma_i
      + h
    \biggr),
  \end{align}
  and set $H_N^{\ne k}(\sigma) \coloneqq H_N(\sigma) - H_N^k(\sigma)$. Denoting by $\sigma^k \in \cS_N$ the configuration that is obtained from $\sigma$ by flipping the spin at site $k \in \{1, \ldots, N\}$, we get $H_N^{\ne k}(\sigma) = H_N^{\ne k}(\sigma^k)$ and $H_N^k(\sigma) = -H_N^k(\sigma^k)$. Since for any $\sigma \in \cS_{N}$ the random variables $H_N^{\ne k}(\sigma)$ and $H_N^k(\sigma)$ are conditionally independent given $\mathcal{G}$, it follows that
  \begin{align}
    \label{eq:EB:H_NvH_N}
    \mean_{\mathcal{G}}\Bigl[
      \me^{\beta(H_N(\sigma) \vee H_N(\sigma^k))} 
    \Bigr]
    &=
    \mean_{\mathcal{G}}\Bigl[ \me^{\beta H_N^{\ne k}(\sigma)} \Bigr]
    \mean_{\mathcal{G}}\Bigl[ 
      \me^{\beta H_N^k(\sigma)} \vee \me^{\beta H_N^k(\sigma^k)} 
    \Bigr].
  \end{align}
  In order to estimate the second term, note that for any two non-negative random variables $X, Y$ with finite second moment, we have
  \begin{align}
    \mean[X \vee Y]
    &=
    \mean[X] \vee \mean[Y]
    + \frac{1}{2}
    \Bigl(
      \mean\bigl[\vert X - Y\vert\bigr] - \bigl\vert \mean[X] - \mean[Y] \bigr\vert
    \Bigr)
    \nonumber\\
    &\leq
    \mean[X] \vee \mean[Y]
    + \frac{1}{2}
    \Bigl(
      \mean\bigl[\vert X - \mean[X]\vert \bigr] + \mean\bigl[\vert Y - \mean[Y]\vert \bigr]
    \Bigr)
    \nonumber\\
    &\leq
    \mean[X] \vee \mean[Y]
    + \frac{1}{2}
    \Bigl( \sqrt{\Var[X]} + \sqrt{\Var[Y]} \Bigr),
  \end{align}
  where we use the triangular inequality and Jensen's inequality. Hence, 
  \begin{align}
    &\mean_{\mathcal{G}}\Bigl[ 
      \me^{\beta H_N^k(\sigma)} \vee \me^{\beta H_N^k(\sigma^k)} 
    \Bigr]
    \nonumber\\
    &\mspace{12mu}\leq
    \mean_{\mathcal{G}}\Bigl[\me^{\beta H_{N}^{k}(\sigma)}\Bigr]
    \vee \mean_{\mathcal{G}}\Bigl[\me^{\beta H_{N}^{k}(\sigma^{k})}\Bigr]
    + \frac{1}{2}
    \biggl(
      \sqrt{\Var_{\mathcal{G}}\Bigl[\me^{\beta H_N^k(\sigma)}\Bigr]}
      +
      \sqrt{\Var_{\mathcal{G}}\Bigl[\me^{\beta H_N^k(\sigma^k)}\Bigr]}
    \biggr).
  \end{align}
  Since,
  \begin{align}
    1+\frac{\Var_{\mathcal{G}}\bigl[\me^{\pm \beta H_{N}^{k}(\sigma)}\bigr]}
    {\mean_{\mathcal{G}}\bigl[\me^{\pm \beta H_{N}^{k}(\sigma)}\bigr]^{2}}
    &\leq
    \frac{\mean_{\mathcal{G}}\Bigl[
        \me^{\pm 2 \beta (H_{N}^{k}(\sigma) - \mean_{\cG}[H_N^k(\sigma)])}
      \Bigr]}{\mean_{\mathcal{G}}\Bigl[
        \me^{\pm \beta (H_{N}^{k}(\sigma) - \mean_{\cG}[H_N^k(\sigma)])}
      \Bigr]^2} 
    \nonumber\\
    &=
    \frac{\exp\biggl(
        \sum_{j=k+1}^N
        \Lambda_{kj}(\tfrac{\mp 2\beta}{N}\sigma_k\sigma_j)
        +
        \sum_{i=1}^{k-1}
        \Lambda_{ik}(\tfrac{\mp 2 \beta}{N}\sigma_i \sigma_k)
      \biggr)}{\exp\biggl(2
        \sum_{j=k+1}^N
        \Lambda_{kj}(\tfrac{\mp \beta}{N}\sigma_k\sigma_j)
        +
        2\sum_{i=1}^{k-1}
        \Lambda_{ik}(\tfrac{\mp \beta}{N}\sigma_i \sigma_k)
      \biggr)},
  \end{align}
  and $\vert \Lambda_{ij}(t)\vert  \leq k_J^2 t^2$, which follows from a Taylor expansion of $\Lambda_{ij}(t)$ up to second order together with the estimate $\vert \Lambda_{ij}''(t)\vert  \leq 2 k_J^2$, we obtain that
  \begin{align}
    \label{eq:EB:H_N^kvH_N^k}
    \mean_{\mathcal{G}}\Bigl[
      \me^{\beta H_N^k(\sigma)} \vee \me^{\beta H_N^k(\sigma^k)} 
    \Bigr]
    \leq
    \biggl(
      \mean_{\mathcal{G}}\Bigl[\me^{\beta H_{N}^{k}(\sigma)}\Bigr]
      \vee \mean_{\mathcal{G}}\Bigl[\me^{\beta H_{N}^{k}(\sigma^{k})}\Bigr]
    \biggr)
    \Bigl(1 + \sqrt{\me^{\frac{6(\beta k_J)^2}{N}}-1}\Bigr).
  \end{align}
  Combining \eqref{eq:EB:H_N^kvH_N^k}, \eqref{eq:EB:H_NvH_N} and \eqref{eq:EB:H}, we see that there exists a $c \equiv c(\beta, k_J)$ such that for all $2\leq N \in \mathbb{N}$,
  \begin{align}
    \mean_{\mathcal{G}}\Bigl[ \me^{\beta(H_N(\sigma) \vee H_N(\sigma^k))} \Bigr]
    \leq
    \me^{\beta(\widetilde{H}_N(\sigma) \vee \widetilde{H}_N(\sigma^k))}\,
    \me^{\alpha_N}\, \Bigl(1 + \frac{c}{\sqrt{N}}\Bigr).
  \end{align}
  This concludes the proof of \eqref{eq:EB:HvH}.
\end{proof}

We are ready to prove the annealed capacity estimates.

\begin{proposition}\label{prop:EB:cap:estimates}
  Let $2 \leq N \in \mathbb{N} $, and let $\mathcal{X}, \mathcal{Y} \subset \cS_N$ be two non-empty and disjoint.
  \begin{enumerate}[(i)]
  \item Then, 
    \begin{align}
      \label{eq:EB:cap:log}
      \Bigl\vert 
        \mean_{\mathcal{G}}\bigl[
          \log\bigl(Z_N \capa_N(\mathcal{X}, \mathcal{Y})\bigr)
        \bigr]
        - \log\bigl(\widetilde{Z}_N \widetilde{\capa}_N(\mathcal{X}, \mathcal{Y})\bigr)
      \Bigr\vert 
      =
      \alpha_N + O\Bigl(\frac{1}{\sqrt{N}}\Bigr).
    \end{align}
  \item For any $q \in [1, \infty)$ there exists a $c_3 \in (0, \infty)$ such that
    \begin{align}
      \label{eq:EB:cap:moments}
      \me^{-\alpha_N} \Bigl(1 - \frac{c_3}{\sqrt{N}}\Bigr)
      &\leq
      \frac{
        \mean_{\mathcal{G}}\Bigl[
          \bigl(Z_N \capa_N(\mathcal{X}, \mathcal{Y})\bigr)^q
        \Bigr]^{1/q}
      }
      {\widetilde{Z}_N \widetilde{\capa}_N(\mathcal{X}, \mathcal{Y})}
      \leq
      \me^{q \alpha_N} \Bigl(1 + \frac{c_3}{\sqrt{N}}\Bigr),
      \\
      \label{eq:EB:1/cap:moments}
      \me^{-\alpha_N} \Bigl(1 - \frac{c_3}{\sqrt{N}}\Bigr)
      &\leq
      \frac{
        \mean_{\mathcal{G}}\Bigl[
          \bigl(Z_N \capa_N(\mathcal{X}, \mathcal{Y})\bigr)^{-q}
        \Bigr]^{1/q}
      }{
        \bigl(
          \widetilde{Z}_N \widetilde{\capa}_N(\mathcal{X}, \mathcal{Y})
        \bigr)^{{-1}}
      }
      \leq
      \me^{q \alpha_N} \Bigl(1 + \frac{c_3}{\sqrt{N}}\Bigr),
    \end{align}
  \end{enumerate}
  where $\alpha_N$ is defined in \eqref{eq:defAlphaN}.
\end{proposition}

\begin{proof}
  Fix $\mathbb{N} \ni N \geq 2$ and consider two non-empty disjoint subsets $\mathcal{X}, \mathcal{Y} \subset \cS_{N}$. Recall from Section~\ref{sec:meth} the definition of the Dirichlet form $\mathcal{E}_{N}(f)$ for functions $f \in \mathcal{H}_{\mathcal{X}, \mathcal{Y}}$ and the Dirichlet form $\mathcal{D}_{N}(\varphi)$ for unit flows $\varphi \in \mathcal{U}_{\mathcal{X}, \mathcal{Y}}$. In view of Lemma~\ref{lem:EB:Delta+HvH}(ii) we have that
  \begin{align}
    \label{eq:EB:forms:estimate}
    \begin{split}
      \mean_{\mathcal{G}}\bigl[Z_{N} \mathcal{E}_{N}(f)\bigr]
      &=
      \widetilde{Z}_{N} \widetilde{\mathcal{E}}_{N}(f)\, 
      \me^{\alpha_{N}} \bigl(1+O(N^{-1/2})\bigr)
      \mspace{72mu} \forall\, f \in \mathcal{H}_{\mathcal{X}, \mathcal{Y}},
      \\
      \mean_{\mathcal{G}}\bigl[Z_{N}^{-1} \mathcal{D}_{N}(\varphi)\bigr]
      &=
      \widetilde{Z}_{N}^{-1} \widetilde{\mathcal{D}}_{N}(\varphi)\, 
      \me^{\alpha_{N}} \bigl(1+O(N^{-1/2})\bigr)
      \mspace{59mu} \forall\, \varphi \in \mathcal{U}_{\mathcal{X}, \mathcal{Y}}.
    \end{split}
  \end{align}
  
  \noindent
  \textit{(i)} The claim in \eqref{eq:EB:cap:log} is an immediate consequence of the Dirichlet principle and the Thomson principle combined with Jensen's inequality.  Indeed, in view of \eqref{eq:EB:forms:estimate} there exists a $c \equiv c(\beta, k_J)$ such that
  \begin{align}
    \mean_{\mathcal{G}}\bigl[
      \log\bigl(Z_{N} \capa_{N}(\mathcal{X}, \mathcal{Y})\bigr)
    \bigr]
    &\leq
    \inf_{f \in \mathcal{H}_{{\mathcal{X}, \mathcal{Y}}}} 
    \log \mean_{\mathcal{G}}\bigl[Z_{N} \mathcal{E}_{N}(f)\bigr]
    \nonumber\\
    &\leq
    \inf_{f \in \mathcal{H}_{{\mathcal{X}, \mathcal{Y}}}} 
    \log\bigl(\widetilde{Z}_{N} \widetilde{\mathcal{E}}_{N}(f)\bigr)
    + \alpha_{N} + c N^{-1/2}
    \nonumber\\
    &=
    \log\bigl(
      \widetilde{Z}_{N}\,\widetilde{\capa}_{N}(\mathcal{X}, \mathcal{Y})
    \bigr)
    + \alpha_{N} + c N^{-1/2}.
    \label{eq:EB:cap:log:ub}
  \end{align}
  Likewise, we obtain that
  \begin{align}
    \mean_{\mathcal{G}}\bigl[
      \log\bigl(Z_{N} \capa_{N}(\mathcal{X}, \mathcal{Y})\bigr)
    \bigr]
    &\geq
    - \inf_{\varphi \in \mathcal{U}_{\mathcal{X}, \mathcal{Y}}}
    \log \mean_{\mathcal{G}}\bigl[Z_{N}^{-1} \mathcal{D}_{N}(\varphi)\bigr]
    \nonumber\\
    &\geq
    - \inf_{\varphi \in \mathcal{U}_{\mathcal{X}, \mathcal{Y}}}
    \log\bigl(\widetilde{Z}_{N}^{-1} \widetilde{\mathcal{D}}_{N}(\varphi)\bigr)
    - \alpha_{N} - c\, N^{-1/2}
    \nonumber\\
    &=
    \log\bigl(
      \widetilde{Z}_{N}\, \widetilde{\capa}_{N}(\mathcal{X}, \mathcal{Y})
    \bigr)
    - \alpha_{N} - c N^{-1/2}.
    \label{eq:EB:cap:log:lb}
  \end{align}
  
  \noindent
  \textit{(ii)} Since the proofs of \eqref{eq:EB:cap:moments} and \eqref{eq:EB:1/cap:moments} are similar, we present the proof for \eqref{eq:EB:1/cap:moments} only. To get the lower bound, note that by Jensen's inequality it is immediate that
  \begin{align}
    \mean_{\mathcal{G}}\Bigl[
      \bigl(Z_{N} \capa_{N}(\mathcal{X}, \mathcal{Y})\bigr)^{-q}
    \Bigr]^{1/q}
    \geq
    \mean_{\mathcal{G}}\Bigl[
      \bigl(Z_{N} \capa_{N}(\mathcal{X}, \mathcal{Y})\bigr)^{-1}
    \Bigr]
    \geq
    \frac{1}{
      \mean_{\mathcal{G}}\bigl[
        Z_{N}\capa_{N}(\mathcal{X}, \mathcal{Y})
      \bigr]
    }.
  \end{align}
  Hence, analogously to \eqref{eq:EB:cap:log:ub}, by applying the Dirichlet principle and \eqref{eq:EB:forms:estimate} we obtain that there exists a $c \equiv c(\beta, k_J)$ such that
  \begin{align}
    \frac{1}
    {\mean_{\mathcal{G}}\bigl[Z_{N}\capa_{N}(\mathcal{X}, \mathcal{Y})\bigr]}
    \geq
    \frac{\me^{-\alpha_N}}{
      \widetilde{Z}_N\, \widetilde{\capa}_N(\mathcal{X}, \mathcal{Y})
    } \bigl(1 - c N^{-1/2}\bigr).
  \end{align}
  To get the upper bound, note that, analogously to \eqref{eq:EB:cap:log:lb}, by the Thomson principle we have that
  \begin{align}
    \mean_{\mathcal{G}}\Bigl[
      \bigl(Z_{N} \capa_{N}(\mathcal{X}, \mathcal{Y})\bigr)^{-q}
    \Bigr]^{1/q}
    \leq
    \inf\nolimits_{\varphi \in \mathcal{U}_{\mathcal{X}, \mathcal{Y}}}
    \mean_{\mathcal{G}}\Bigl[
      \bigl( Z_{N}^{-1} \mathcal{D}_{N}(\varphi) \bigr)^{q}
    \Bigr]^{1/q}.
  \end{align}
  By Minkowski's inequality and an application of \eqref{eq:EB:HvH} with $\beta$ replaced by $\beta q$, we find that for any $q\in [1,\infty)$ there exists a $c'\equiv c'(q, \beta, k_J)$ such that, for all $\varphi \in \mathcal{U}_{\mathcal{X}, \mathcal{Y}}$,
  \begin{align}
    \mean_{\mathcal{G}}\Bigl[\bigl(Z_N^{-1} \mathcal{D}_N(\varphi)\bigr)^q\Bigr]^{1/q}
    \leq
    \widetilde{Z}_N^{-1} \widetilde{\mathcal{D}}_N(\varphi)\, \me^{q \alpha_N}\,
    \bigl(1 + c' N^{-1/2}\bigr).
  \end{align}
  Therefore, again applying the Thomson principle, we obtain
  \begin{align}
    \mean_{\mathcal{G}}\Bigl[
      \bigl(Z_{N} \capa_{N}(\mathcal{X}, \mathcal{Y})\bigr)^{-q}
    \Bigr]^{1/q}
    \leq
    \frac{\me^{q \alpha_N}}
    {\widetilde{Z}_N\, \widetilde{\capa}_N(\mathcal{X}, \mathcal{Y})}
    \bigl(1 + c' N^{-1/2}\bigr),
  \end{align} 
  and by setting $c_3 \coloneqq c \vee c'$ we conclude the proof.
\end{proof}


\section{Equilibrium potential estimates} \label{sec:harm}

This section contains all our results concerning the $\ell^1(\mu_N)$-norm of the equilibrium potential $h_{\mathcal{A}_N, \mathcal{B}_N}^N$, which we call the \emph{harmonic sum}. Before proving concentration estimates in Section~\ref{sec:h_conc} and annealed estimates in Section~\ref{sec:h_ann}, we provide some preliminary estimate in Section~\ref{sec:pre_harm}. We emphasise that throughout this section, contrary to Section~\ref{sec:capa}, metastability plays an essential role.


\subsection{Preliminary estimates} \label{sec:pre_harm}
As mentioned above, for estimates on the harmonic sum we restrict to the event $\Xi(a_N)$, which is defined in \eqref{eq:def:Xi(N)} and used in Section~\ref{sec:meta_proof} to prove Theorem~\ref{thm:meta}. This event has high probability for suitably chosen sequences $(a_N)_{N \in \mathbb{N}}$ (recall Lemma~\ref{lemma:conc:H_N}). We use two facts: on $\Xi(a_N)$ we can control the quenched Gibbs measure $\mu_N$ in terms of the annealed Gibbs measure $\tilde{\mu}_N$ (recall \eqref{eq:comparison:mu}), and the harmonic sum localises on the metastable valley of $\cA_N$. We state and prove the last result in Proposition~\ref{prop:harm:estimate}.

\begin{notation}
  For any two random variables $X,Y$ depending on $N\in \mathbb{N}$, writing \newquote{$X=Y$ holds  $\prob$-a.s.\ on the event $\Otimeta(N)$} means that
  \begin{align}
    \ind_{\Otimeta(N)} X= \ind_{\Otimeta(N)} Y, \qquad \prob\text{-a.s.}.
  \end{align}
\end{notation}

We stress that all formulas in Section~\ref{sec:harm} hold $\prob$-a.s. Moreover, formulas involving the quantities $\cS_{j,N}, \mathcal{M}_{j,N}, \mathcal{A}_N, \mathcal{B}_N$, for fixed $j \in \{1, \dots, K\}$ and $N \in \mathbb{N}$, hold $\prob$-a.s. on the event $\Otimeta(N)$ unless differently specified, because those quantities are not defined in $\Otimeta(N)^c$.

Furthermore, recall that the sets $\mathcal{A}_N$ and $\mathcal{B}_N$ are defined, in terms of the fixed index $i \in \{2, \ldots, K\}$, in \eqref{def:A_N:B_N}.

\begin{remark} 
  \label{rem:noInd}
  By $\cG$-measurability of $\Otimeta(N)$ (see Remark~\ref{rem:Gmeas}) we are allowed to compute expectations and probabilities conditioned to $\cG$ on the event $\Otimeta(N)$.
\end{remark}

\begin{proposition}
  \label{prop:harm:estimate}
  Suppose that $(a_N)_{N \in \mathbb{N}}$ is a non-negative sequence that is sublinear in $N$. Then there exists a $C \in (0, k_1 \wedge k_2)$ such that $\prob$-a.s.\ on the event $\Xi(a_N) \cap \Otimeta(N)$, for $N$ sufficiently large depending on $(a_N)_{N \in \mathbb{N}}, \beta, k_1, k_2$, 
  \begin{align}
    \label{eq:harm:estimate}
    \norm{h_{\mathcal{A}_N, \mathcal{B}_N}^N}_{\mu_N}
    = \mu_N[\cS_{i,N}] \bigl[1 + O(\me^{-C N})\bigr].
  \end{align}
\end{proposition}

\noindent

\begin{remark}
  Proposition~\ref{prop:harm:estimate} holds true also for $a_N=c N$, with $c>0$, in case $c \beta$ is sufficiently small compared to $k_1$ and $k_2$. We do not use this result.
\end{remark}

\begin{remark}
  Although Proposition \ref{prop:harm:estimate} is similar to \cite[Theorem~1.7]{SS19}, it is not an immediate consequence of the latter. Indeed, in \eqref{eq:harm:estimate} both $\norm{h_{\mathcal{A}_N, \mathcal{B}_N}^N}_{\mu_N}$ and $\mu_N$ refer to the \emph{quenched} Markov chain $(\Sigma_N(t))_{t \in \mathbb{N}_0}$, but $\cS_{i,N}$  is a set of the metastable partition of the \emph{annealed} Markov chain $(\widetilde{\Sigma}_N(t))_{t \in \mathbb{N}_0}$, while in \cite[Theorem~1.7]{SS19} all quantities refer to the same process. We made this modification for two reasons. First of all, we are not able to prove \cite[Theorem~1.7]{SS19} for the \emph{quenched} Markov chain because we cannot prove the non-degeneracy assumption needed therein. Second, even if we were able to prove it, it would not be useful later on as we do not have estimates on the measure of the metastable partition of $(\Sigma_N(t))_{t \in \mathbb{N}_0}$.  However, as we shall see later, \eqref{eq:non-degeneracy} and \eqref{eq:comparison:mu} allow us both to prove \eqref{eq:harm:estimate} and to use it later having estimates on its right hand side.
\end{remark}

\begin{remark}
  Note that in the following proof of Proposition~\ref{prop:harm:estimate} we use the metastability of $(\widetilde{\Sigma}_N(t))_{t \in \mathbb{N}_0}$ and do not use the metastability of $(\Sigma_N(t))_{t \in \mathbb{N}_0}$.
\end{remark}

\begin{proof}[Proof of Proposition~{\upshape\ref{prop:harm:estimate}}]
  The proof is inspired by that of \cite[Lemma~3.3]{SS19} and consists of two steps.
  
  \noindent  
  \textit{Step 1.} 
  Let $(a_N)_{N \in \mathbb{N}}$ be any non-negative real sequence and fix $N \in \mathbb{N}$. We start by showing that the following is true $\prob$-a.s.\ on the event $\Xi(a_N) \cap \Otimeta(N)$, for any two $j \neq k \in \{1, \dots, K\}$ and any $\varepsilon \in (0,1]$,
  \begin{align}
    \label{eq:harm:claim:1}
    \sum_{\sigma \in \mathcal{S}_{k,N}} 
    \frac{\mu_N(\sigma)}{\mu_N[\mathcal{S}_{k,N}]}\, 
    h_{\mathcal{M}_{j,N}, \mathcal{M}_{k,N}}^N(\sigma)
    \leq
    \varepsilon + \tilde{\rho}_N\, \me^{4 \beta a_N} \log(1/\varepsilon) 
    \min\biggl\{
      1, \frac{\mu_N[\mathcal{S}_{j,N}]}{\mu_N[\mathcal{S}_{k,N}]}
    \biggr\}.
  \end{align}
  Recall that, for $\ell \in \{1, \ldots, K\}$, $\mathcal{V}_{\ell ,N}$ and $\cS_{\ell ,N}$ are, respectively, the local valley around the metastable set $\mathcal{M}_{\ell ,N}$ and a set of the metastable partition of the \emph{annealed model}. Moreover, $\mathcal{M}_N=\bigcup_{\ell=1}^K \mathcal{M}_{\ell,N}$. By applying \eqref{eq:P_hit_cond} and \eqref{eq:comparison:exit_prob}, we have that, for any $\mathcal{X} \subset \mathcal{V}_{k,N} \setminus \mathcal{M}_{k,N}$, $\prob$-a.s.\ on the event  $\Xi(a_N) \cap \Otimeta(N)$,
  \begin{align}
    \mu_N&\bigl[\mathcal{X}\bigr]
    \overset{\mspace{-15mu}\eqref{eq:comparison:exit_prob}\mspace{-15mu}}{\leq}
    \me^{2\beta a_N}
    \frac{\capa_N(\mathcal{X}, \mathcal{M}_{k,N})}
    {
      \tProb_{\mu_N \vert  \mathcal{X}}^N\bigl[
        \tilde\tau_{\mathcal{M}_{k,N}}^N < \tilde\tau_{\mathcal{X}}^N
      \bigr]
    }
    \nonumber\\
    &\leq
    \me^{2 \beta a_N}
    \tilde{\rho}_N
    \biggl(
      \max_{\ell \in \{1, \ldots, K\}} 
      \tProb_{\mu_N \vert  \mathcal{M}_{\ell,N}}^N\Bigl[
        \tilde\tau_{\mathcal{M}_N \setminus \mathcal{M}_{\ell,N}}^N 
        < \tilde\tau_{\mathcal{M}_{\ell,N}}^N
      \Bigr]
    \biggr)^{-1}
    \capa_N(\mathcal{X}, \mathcal{M}_{k,N})
    \nonumber\\ 
    &\overset{\mspace{-15mu}\eqref{eq:comparison:exit_prob}\mspace{-15mu}}{\leq}
    \me^{4 \beta a_N}
    \tilde{\rho}_N
    \biggl(
      \max_{\ell \in \{1, \ldots, K\}} 
      \Prob_{\mu_N \vert  \mathcal{M}_{\ell,N}}^N\Bigl[
        \tau_{\mathcal{M}_N \setminus \mathcal{M}_{\ell,N}}^N 
        < \tau_{\mathcal{M}_{\ell,N}}^N
      \Bigr]
    \biggr)^{-1}
    \capa_N(\mathcal{X}, \mathcal{M}_{k,N}),
    \label{eq:muX}
  \end{align}
  where in the second inequality we used that we are in $\Otimeta(N)$ and $\mathcal{X} \subset \mathcal{V}_{k,N} \setminus \mathcal{M}_{k,N}$ to apply \cite[Lemma~3.1]{SS19}. Moreover, using \eqref{eq:P_hit_cond} and monotonicity of capacities, we get
  \begin{align}
    &\max_{\ell \in \{1, \ldots, K\}}
    \Prob_{\mu_N \vert  \mathcal{M}_{\ell,N}}^N\Bigl[
      \tau_{\mathcal{M}_N \setminus \mathcal{M}_{\ell,N}}^N 
      < \tau_{\mathcal{M}_{\ell,N}}^N
    \Bigr]
    \nonumber\\ 
    &\mspace{36mu}\geq
    \max\biggl\{
      \frac{\capa_N(\mathcal{M}_{j,N}, \mathcal{M}_{k,N})}
      {\mu_N\bigl[\mathcal{M}_{k,N}\bigr]},
      \frac{\capa_N(\mathcal{M}_{k,N}, \mathcal{M}_{j,N})}
      {\mu_N\bigl[\mathcal{M}_{j,N}\bigr]}
    \biggr\}.
    \label{eq:maxPMj}
  \end{align}
  Next, for $t \in (0,1]$ we write $\mathcal{X}_N(t) \coloneqq \{\sigma \in \cS_N : h_{\mathcal{M}_{j,N}, \mathcal{M}_{k,N}}^N(\sigma) \geq t\}$ to denote the super level-sets of $h_{\mathcal{M}_{j,N}, \mathcal{M}_{k,N}}^N$.  Note that, for any $t \in (0,1]$, $\mathcal{X}_N(t) \cap \mathcal{M}_{k,N} = \emptyset$ and $\mathcal{M}_{j,N} \subseteq \mathcal{X}_N(t)$. Using reversibility, \eqref{eq:def_capa} and the definition of $\mathcal{X}_N(t)$ we have
  \begin{align}
    t \capa_N(\mathcal{X}_N(t), \mathcal{M}_{k,N})
    &\leq
    \sum_{\sigma \in \mathcal{X}_N(t)} \mu_N(\sigma)
    \bigl(-L_N h_{\mathcal{X}_N(t), \mathcal{M}_{k,N}}^{N}\bigr)(\sigma)\,
    h_{\mathcal{M}_{j,N}, \mathcal{M}_{k,N}}^{N}(\sigma)
    \nonumber\\ 
    &=
    \scpr{-L_N h_{\mathcal{X}_N(t), \mathcal{M}_{k,N}}^N, h_{\mathcal{M}_{j,N}, \mathcal{M}_{k,N}}^N}_{\mu_N}
    \nonumber\\[1.5
    ex] 
    &=
    \scpr{h_{\mathcal{X}_N(t), \mathcal{M}_{k,N}}^N, -L_N h_{\mathcal{M}_{j,N}, \mathcal{M}_{k,N}}^N}_{\mu_N}
    =
    \capa_N(\mathcal{M}_{j, N}, \mathcal{M}_{k,N}).
    \label{eq:tcap}
  \end{align}
  By expressing the expected value of a non-negative random variable in terms of the integral of the tail of its distribution, we obtain, for any $\varepsilon \in (0,1]$,
  \begin{align}
    \sum_{\sigma \in \cS_{k,N}} 
    \frac{\mu_N(\sigma)}{\mu_N[\mathcal{S}_{k,N}]} 
    h_{\mathcal{M}_{j,N}, \mathcal{M}_{k,N}}^N(\sigma)
    \leq 
    \varepsilon
    + 
    \int_{\varepsilon}^1 
    \frac{\mu_N\bigl[\mathcal{X}_N(t) \cap \cS_{k,N}\bigr]}
    {\mu_N\bigl[\cS_{k,N}\bigr]}
    \md t.
  \end{align}
  Using \eqref{eq:muX} with $\mathcal{X} = \mathcal{X}_N(t)\cap  \cS_{k,N}$, together with \eqref{eq:maxPMj}, the symmetry and monotonicity of capacities and \eqref{eq:tcap}, we obtain, for $t \in [\varepsilon, 1]$,
  \begin{align}
    \mu_N\bigl[\mathcal{X}_N(t) \cap \cS_{k,N}\bigr] 
    &\leq 
    \tilde{\rho}_N\, \me^{4 \beta a_N}
    \frac{
      \min\bigl\{
        \mu_N\bigl[\mathcal{M}_{k,N}\bigr],
        \mu_N\bigl[\mathcal{M}_{j,N}\bigr]
      \bigr\}
    }
    {
      \capa_N(\mathcal{M}_{j,N}, \mathcal{M}_{k,N})
    }
    \capa_N(\mathcal{X}_N(t), \mathcal{M}_{k,N})
    \nonumber\\ 
    &\leq  
    \tilde{\rho}_N\, \me^{4 \beta a_N}
    \min\bigl\{
      \mu_N\bigl[\mathcal{M}_{k,N}\bigr],
      \mu_N\bigl[\mathcal{M}_{j,N}\bigr]
    \bigr\}
    \frac{1}{t}.
  \end{align}
  Therefore, recalling that $\mathcal{M}_{\ell,N}\subseteq \mathcal{S}_{\ell,N}$, $\ell \in \{1, \dots, K\}$, we obtain
  \begin{align}
    \sum_{\sigma \in \cS_{k,N}} 
    \frac{\mu_N(\sigma)}{\mu_N[\mathcal{S}_{k,N}]}\, 
    h_{\mathcal{M}_{j,N}, \mathcal{M}_{k,N}}^N(\sigma)
    \leq
    \varepsilon
    +
    \tilde{\rho}_N\, \me^{4\beta a_N}
    \min\biggl\{
      \frac{\mu_N\bigl[\cS_{j,N}\bigr]}{\mu_N\bigl[\cS_{k,N}\bigr]}, 1
    \biggr\}
    \int_{\varepsilon}^1 \frac{1}{t} \md t,
  \end{align} 
  which completes the proof of \eqref{eq:harm:claim:1}.
  
  \noindent 
  \textit{Step 2.} In view of \eqref{eq:harm:claim:1}, the proof of \eqref{eq:harm:estimate} runs along the same lines as the proof of \cite[Theorem~1.7]{SS19}.  For the reader's convenience we provide the details here. Let $\mathcal{A}_N, \mathcal{B}_N$ be defined as in \eqref{def:A_N:B_N}. In particular, recall that $\mathcal{A}_N = \mathcal{M}_{i,N}$. Then 
  \begin{align}
    \label{eq:hABmu_dec}
    \norm{h_{\mathcal{A}_N, \mathcal{B}_N}^N}_{\mu_N}
    =
    \mu_N\bigl[\cS_{i,N}\bigr]
    \Biggl(
      \norm{h_{\mathcal{A}_{N}, \mathcal{B}_N}^N}_{\mu_N \vert  \cS_{i,N}}
      +
      \sum_{j \ne i}
      \frac{\mu_N\bigl[\cS_{j,N}\bigr]}{\mu_N\bigl[\cS_{i,N}\bigr]}
      \norm{h_{\mathcal{A}_{N}, \mathcal{B}_N}^N}_{\mu_N \vert  \cS_{j,N}}
    \Biggr).
  \end{align}
  In order to prove a lower bound, we neglect the last term in the bracket in \eqref{eq:hABmu_dec}, while the first term is bounded from below by
  \begin{align}
    \norm{h_{\mathcal{A}_N, \mathcal{B}_N}^N}_{\mu_N \vert  \cS_{i,N}}
    =
    1 - \norm{h_{\mathcal{B}_N, \mathcal{A}_N}^N}_{\mu_N \vert  \cS_{i,N}}
    \geq
    1 - \sum_{j=1}^{i-1} 
    \norm{h_{\mathcal{M}_{j,N}, \mathcal{M}_{i,N}}^N}_{\mu_N \vert  \cS_{i,N}},
  \end{align}
  where we used that, for all $\sigma \in \cS_N \setminus (\mathcal{A}_N \cup \mathcal{B}_N)$,
  \begin{align}
    \label{eq:splitBA}
    \begin{split}
      h_{\mathcal{B}_N, \mathcal{A}_N}^N(\sigma)
      &=
      \Prob_{\sigma}^N\Bigl[
        \tau_{\bigcup_{j=1}^{i-1} \mathcal{M}_{j,N}}^N
        < \tau_{\mathcal{M}_{i,N}}^N
      \Bigr]
      \\ 
      &\leq
      \sum_{j=1}^{i-1}
      \Prob_{\sigma}^N\Bigl[
        \tau_{\mathcal{M}_{j,N}}^N <  \tau_{\mathcal{M}_{i,N}}^N
      \Bigr]
      =
      \sum_{j=1}^{i-1} h_{\mathcal{M}_{j,N}, \mathcal{M}_{i,N}}^N(\sigma).
    \end{split}
  \end{align}
  By applying \eqref{eq:harm:claim:1} with $\varepsilon = \me^{-k_1 N}$, recalling $\tilde{\rho}_N=\me^{-k_1 N}$ (see \eqref{eq:def_tirho}), we obtain that $\prob$-a.s.\ on the event $\Xi(a_N) \cap \Otimeta(N)$,
  \begin{align}
    \label{eq:1-K}
    \norm{h_{\mathcal{A}_N, \mathcal{B}_N}^N}_{\mu_N \vert  \cS_{i,N}}
    \geq
    1 - K \me^{-k_1 N}
    \bigl(1 + \me^{4 \beta a_N} \log(1/\me^{-k_1 N})\bigr).
  \end{align} 
  Hence,  we get
  \begin{align}
    \label{eq:harm:sum:lb}
    \norm{h_{\mathcal{A}_N, \mathcal{B}_N}^N}_{\mu_N}
    \geq
    \mu_N\bigl[\cS_{i,N}\bigr]
    \Bigl(
      1 - K N \me^{-k_1 N +\beta a_N} \bigl(\me^{-\log N -\beta a_N} + k_1\bigr)
    \Bigr).
  \end{align}
  To get the upper bound, we exploit the fact that \eqref{eq:non-degeneracy}, together with \eqref{eq:comparison:mu}, implies that $\mu_N\bigl[\cS_{j,N}\bigr] / \mu_N\bigl[\cS_{i,N}\bigr] \leq \me^{-k_2N}\, \me^{2\beta a_N}$ for all $j \in \{i+1, \ldots, K\}$. Hence
  \begin{align}
    &\sum_{j \ne i}
    \frac{\mu_N\bigl[\cS_{j,N}\bigr]}{\mu_N\bigl[\cS_{i,N}\bigr]}
    \norm{h_{\mathcal{A}_N, \mathcal{B}_N}^N}_{\mu_N \vert  \cS_{j,N}} \\ \nonumber
    &\qquad \leq
    K \me^{-k_2N}\, \me^{2\beta a_N}
    +
    \sum_{j=1}^{i-1}
    \frac{\mu_N\bigl[\cS_{j,N}\bigr]}{\mu_N\bigl[\cS_{i,N}\bigr]}
    \norm{h_{\mathcal{M}_{i,N}, \mathcal{M}_{j,N}}^N}_{\mu_N \vert  \cS_{j,N}},
  \end{align}
  where we used that, for $j \in \{1, \ldots, i-1\}$ and $\sigma \in \cS_N \setminus \bigcup_{\ell=1}^i \mathcal{M}_{\ell,N}$, 
  \begin{align}
    h_{\mathcal{A}_N, \mathcal{B}_N}^N (\sigma) 
    &= 
    \Prob_{\sigma}^N\Bigl[  
      \tau_{\mathcal{M}_{i,N}} ^N  
      < \tau_{\bigcup_{\ell=1}^{i-1} \mathcal{M}_{\ell,N}}^N   
    \Bigr] 
    \nonumber \\
    &\leq 
    \Prob_{\sigma}^N\Bigl[
      \tau_{\mathcal{M}_{i,N}} ^N  < \tau_{ \mathcal{M}_{j,N}}^N   
    \Bigr] 
    = 
    h_{\mathcal{M}_{i,N}, \mathcal{M}_{j,N}}^N (\sigma).
  \end{align}
  Thus, applying \eqref{eq:harm:claim:1} with $\varepsilon = \me^{-k_1 N} \min_{\ell \in \{1, \ldots, i-1\}} \mu_N\bigl[\cS_{i,N}\bigr]/\mu_N\bigl[\cS_{\ell,N}\bigr]$, we get that, $\prob$-a.s.\ on the event $\Xi(a_N) \cap \Otimeta(N)$,
  \begin{align}
    \begin{split}
      \frac{\mu_N\bigl[\cS_{j,N}\bigr]}{\mu_N\bigl[\cS_{i,N}\bigr]}&
      \norm{h_{\mathcal{M}_{i,N}, \mathcal{M}_{j,N}}^N}_{\mu_N \vert  \cS_{j,N}}
      \\
      &\leq
      \me^{-k_1 N} 
      - 
      \tilde{\rho}_N\, \me^{4 \beta a_N} 
      \log\biggl(
        \me^{-k_1 N} \min_{\ell \in\{1, \ldots, i-1\}} 
        \frac{\mu_N\bigl[\cS_{i,N}\bigr]}{\mu_N\bigl[\cS_{\ell,N}\bigr]}
      \biggr)
    \end{split}
  \end{align}
  for $j \in \{1, \ldots, i-1\}$. Since $\mu_N\bigl[\cS_{i,N}\bigr] / \mu_N\bigl[\cS_{j,N}\bigr] \geq \me^{-\beta (k_J +h)N}$ for $j \in \{1, \ldots, K\}$ and $\norm{h_{\mathcal{A}_N, \mathcal{B}_N}^N}_{\mu_N \vert  \cS_{i,N}} \leq 1$, we can use \eqref{eq:hABmu_dec} to conclude that
  \begin{align}
    &\norm{h_{\mathcal{A}_N, \mathcal{B}_N}^N}_{\mu_N}
    \nonumber\\
    &\leq
    \mu_N\bigl[\cS_{i,N}\bigr]
    \Bigl(
      1 + K \me^{-k_2N +2\beta a_N}
      + K \bigl(\me^{-k_1 N} + (k_1 + \beta(k_J + h)) N \me^{-k_1 N + 4 \beta a_N} \bigr)
    \Bigr).
    \label{eq:harm:sum:ub}
  \end{align}
  Let 
  \begin{align}
    \bar{N} \coloneqq \min\{N \in \mathbb{N}\colon -k_2N+2\beta a_N< 0 \text{ and } - k_1 N +4\beta a_N+\log N<0\}.
  \end{align}
  Note that $\bar{N}$ depends on $(a_N)_{N \in \mathbb{N}}, \beta, k_1,k_2$ and is deterministic. The minimum exists because $\beta, k_1,k_2$ are fixed and $a_N$ is taken sublinear in $N$.
  By combining \eqref{eq:harm:sum:lb} and \eqref{eq:harm:sum:ub}, the assertion follows for all $N\geq \bar{N}$.
\end{proof}

\begin{corollary}
  \label{cor:mean:harm:estimate}
  There exists a $C \in (0, k_1 \wedge k_2)$ such that, for $N$ sufficiently large depending on $\beta, k_1, k_2, k_J$, $\prob$-a.s.\ on the event $\Otimeta(N)$,  
  \begin{align}
    \label{eq:mean:harm:estimate}
    \mean_{\mathcal{G}}\Bigl[
      \log\bigl(Z_N \norm{h_{\mathcal{A}_N, \mathcal{B}_N}^N}_{\mu_N}\bigr)
    \Bigr]
    =
    \mean_{\mathcal{G}}\Bigl[
      \log\bigl(Z_N \mu_N[\cS_{i,N}]\bigr)
    \Bigr]
    + O\bigl(\me^{-C N}\bigr).
  \end{align}
\end{corollary}

\begin{proof}
  First observe that $\prob$-a.s.\ 
  \begin{align}
    \me^{-\beta (k_J +h) N} \leq Z_N \norm{h_{\mathcal{A}_N, \mathcal{B}_N}^N}_{\mu_N} \leq 2^N \me^{\beta (k_J +h) N}
  \end{align}
  for $N \in \mathbb{N} \setminus\{1\}$.  Moreover, let $a_N = k_J\sqrt{N k_1+ (N+1)\log 2}$. In view of Proposition~\ref{prop:harm:estimate}, we know that there exist a $C \in (0, k_1 \wedge k_2)$ and a $c' \in (0, \infty)$ such that, for all $N$ sufficiently large depending on $\beta, k_1, k_2, k_J$, 
  \begin{align}
    &\mean_{\mathcal{G}}
    \Bigl[
      \log\bigl(Z_N \norm{h_{\mathcal{A}_N, \mathcal{B}_N}^N}_{\mu_N}\bigr)
    \Bigr]
    \nonumber\\ 
    &\leq
    \mean_{\mathcal{G}}\Bigl[
      \log\bigl(Z_N \norm{h_{\mathcal{A}_N, \mathcal{B}_N}^N}_{\mu_N}\bigr)
      \ind_{\Xi(a_N) 
      }
    \Bigr]
    +
    (\beta (k_J +h) + \log 2) N \prob_{\mathcal{G}}\bigl[\Xi(a_N)^c\bigr]
    \nonumber\\ 
    &\leq
    \mean_{\mathcal{G}}\Bigl[
      \log\bigl(Z_N \mu_N\bigl[\cS_{i,N}\bigr]\bigr)
    \Bigr]
    +
    \log\bigl(1 + c'(\me^{-C N})\bigr)
    +
    (\beta (k_J +h) + \log 2) N\, \me^{-k_1 N},
  \end{align}
  where we used \eqref{eq:tail:estimate}, which is implied by Lemma~\ref{lemma:conc:H_N} and our choice of $a_N$. Likewise, by using additionally that $Z_N \mu_N[\cS_{i,N}] \leq 2^N \me^{\beta (k_J +h) N}$, we obtain that
  \begin{align}
    &\mean_{\mathcal{G}}
    \Bigl[
      \log\bigl(Z_N \norm{h_{\mathcal{A}_N, \mathcal{B}_N}^N}_{\mu_N}\bigr)
    \Bigr]
    \nonumber\\ 
    &\geq
    \mean_{\mathcal{G}}\Bigl[
      \log\bigl(Z_N \norm{h_{\mathcal{A}_N, \mathcal{B}_N}^N}_{\mu_N}\bigr)
      \ind_{\Xi(a_N)
      }
    \Bigr]
    -
    \beta (k_J +h) N \prob_{\mathcal{G}}\bigl[\Xi(a_N)^c\bigr]
    \nonumber\\ 
    &\geq 
    \mean_{\mathcal{G}}\Bigl[
      \log\bigl(Z_N \mu_N\bigl[\cS_{i,N}\bigr]\bigr)
    \Bigr]
    +
    \log\bigl(1 - c'(\me^{-C N})\bigr)
    -
    (2\beta (k_J +h) + \log 2) N \me^{-k_1 N}.
  \end{align}
  Since $C < k_1$, this concludes the proof. 
\end{proof}  


\subsection{Concentration inequality} \label{sec:h_conc}

\begin{proposition}
  \label{pro:conc_Zh}
  There exist $C \in (0, k_1 \wedge k_2)$ and $c_4 \in (0, \infty)$ such that, for all $N$ sufficiently large depending on $\beta, k_1, k_2, k_J$, and all $t \in \mathbb{N}_0$,
  $\prob$-a.s.\ on the event $\Otimeta(N)$,
  \begin{align}
    \prob_{\mathcal{G}} \Bigl[
      \bigl\vert 
        \log \bigl(Z_N \norm{h_{\mathcal{A}_N,\mathcal{B}_N}^N}_{\mu_N}\bigr) 
        - 
        \mean_{\mathcal{G}}\bigl[ 
          \log \bigl(Z_N \norm{h_{\mathcal{A}_N,\mathcal{B}_N}^N}_{\mu_N} \bigr) 
        \bigr]
      \bigr\vert  > t
    \Bigr] 
    \leq 
    2 \me^{-\bigl(\frac{t-c_N}{\beta k_J}\bigr)^2} + \me^{-k_1 N}\mspace{-4mu},
    \label{eq:conc_Zh}
  \end{align}
  where $c_N \coloneqq c_4\, \me^{-C N}$.
\end{proposition}

\begin{proof}
  Let $a_N = k_J\sqrt{N k_1+ (N+1)\log 2}$. In view of Proposition~\ref{prop:harm:estimate} and Corollary~\ref{cor:mean:harm:estimate}, there exist a $C \in (0, k_1 \wedge k_2)$ and $c_4 \in (0, \infty)$ such that, for $N$ sufficiently large depending on $\beta, k_1, k_2, k_J$, $\prob$-a.s on the event $\Xi(a_N)\cap \Otimeta(N)$,
  \begin{align}
    \bigl\vert  
      \log\bigl(Z_N \norm{h_{\mathcal{A}_N, \mathcal{B}_N}^N}_{\mu_N}\bigr)
      -
      \log\bigl(Z_N \mu_N[\cS_{i,N}]\bigr)
    \bigr\vert 
    \leq
    \frac{c_4}{2} \, \me^{-C N}
  \end{align}
  and
  \begin{align}
    \bigl\vert 
      \mean_{\mathcal{G}}\bigl[
        \log\bigl(Z_N \norm{h_{\mathcal{A}_N, \mathcal{B}_N}^N}_{\mu_N}\bigr)
      \bigr]
      -
      \mean_{\mathcal{G}}\bigl[
        \log\bigl(Z_N \mu_N[\cS_{i,N}]\bigr)
      \bigr]
    \bigr\vert 
    \leq
    \frac{c_4}{2} \, \me^{-C N}.
  \end{align}
  Hence, by setting $c_N \coloneqq c_4\, \me^{-C N}$, we obtain that
  \begin{align}
    &\prob_{\mathcal{G}}\Bigl[
      \bigl\vert  
        \log\bigl(Z_N \norm{h_{\mathcal{A}_N, \mathcal{B}_N}^N}_{\mu_N}\bigr)
        -
        \mean_{\mathcal{G}}\bigl[
          \log\bigl(Z_N \norm{h_{\mathcal{A}_N, \mathcal{B}_N}^N}_{\mu_N}\bigr)
        \bigr]
      \bigr\vert  > t
    \Bigr]
    \nonumber\\ 
    &\leq
    \prob_{\mathcal{G}}\Bigl[
      \bigl\vert  
        \log\bigl(Z_N \norm{h_{\mathcal{A}_N, \mathcal{B}_N}^N}_{\mu_N}\bigr)
        -
        \mean_{\mathcal{G}}\bigl[
          \log\bigl(Z_N \norm{h_{\mathcal{A}_N, \mathcal{B}_N}^N}_{\mu_N}\bigr)
        \bigr]
      \bigr\vert  > t,
      \Xi(a_N)
    \Bigr]
    \nonumber\\ 
    &\quad 
    +
    \prob_{\mathcal{G}}\bigl[\Xi(a_N)^c\bigr]
    \nonumber\\ 
    &\leq
    \prob_{\mathcal{G}}\Bigl[
      \bigl\vert  
        \log\bigl(Z_N \mu_N[\cS_{i,N}]\bigr)
        -
        \mean_{\mathcal{G}}\bigl[
          \log\bigl(Z_N \mu_N[\cS_{i,N}]\bigr)
        \bigr]
      \bigr\vert  > t - c_N
    \Bigr]
    +
    \me^{-k_1 N},
    \label{eq:harm:conc:estimate}
  \end{align}
  where, as above, we used \eqref{eq:tail:estimate}, which is implied by Lemma~\ref{lemma:conc:H_N} and our choice of $a_N$.  In order to bound the first term on the right-hand side of \eqref{eq:harm:conc:estimate}, recall that the triangular array $(J_{ij})_{1 \leq i < j < \infty}$ is assumed to be conditionally independent given $\mathcal{G}$. Moreover, in view of \eqref{eq:H_N:JvsJ'}, for any $2 \leq N \in \mathbb{N}$ it is immediate that the mapping
  \begin{align}
    (J_{ij})_{1 \leq i < j \leq N}
    \longmapsto
    \bar{F}_N\bigl((J_{ij})_{1 \leq i < j \leq N}\bigr)
    \coloneqq 
    \log\bigl(Z_N \mu_N[\cS_{i,N}]\bigr)
  \end{align}
  satisfies the estimate
  \begin{align}
    \bigl\vert 
      \bar{F}_N\bigl((J_{ij})_{1 \leq i < j \leq N}\bigr)
      -
      \bar{F}_N\bigl((J_{ij}')_{1 \leq i < j \leq N}\bigr)
    \bigr\vert 
    \leq
    \frac{2\beta k_J}{N},
  \end{align}
  where $J_{ij}' \coloneqq J_{ij}$ for all $1 \leq i < j \leq N$ such that $(i,j) \ne (k,l)$ and $J_{kl}'$ is a conditionally independent copy of $(J_{ij})_{1 \leq i < j \leq N}$ given $\mathcal{G}$, for any $1 \leq k < l \leq N$. Hence, by applying McDiarmid concentration inequality in the version of Proposition~\ref{prop:McDiarmid}, we get that 
  \begin{align}
    \label{eq:mu:conc}
    \prob_{\mathcal{G}}\Bigl[
      \bigl\vert  
        \log\bigl(Z_N \mu_N[\cS_{i,N}]\bigr)
        -
        \mean_{\mathcal{G}}\bigl[
          \log\bigl(Z_N \mu_N[\cS_{i,N}]\bigr)
        \bigr]
      \bigr\vert  > t
    \Bigr]
    \leq
    2\, \me^{- t^2/(\beta k_J)^2}.
  \end{align} 
  Combining \eqref{eq:harm:conc:estimate} and \eqref{eq:mu:conc}, we get the assertion in \eqref{eq:conc_Zh}.
\end{proof}


\subsection{Annealed estimate} \label{sec:h_ann}

\begin{proposition}
  \label{prop:harm:annealed}
  Let $\alpha_N$ be as defined in \eqref{eq:defAlphaN}. Then the following hold:
  \begin{enumerate}[(i)]
  \item There exists a $c_5 \in (0, \infty)$ such that, for $N$ sufficiently large depending on $\beta, k_1, k_2, k_J$, $\prob$-a.s.\ on the event $\Otimeta(N)$,
    \begin{align}
      \label{eq:log:harm:annealed}
      - \frac{c_5}{N}
      \leq
      \mean_{\mathcal{G}}\Bigl[
        \log\bigl(Z_N \norm{h_{\mathcal{A}_N,\mathcal{B}_N}^N}_{\mu_N}\bigr)
      \Bigr]
      -
      \log\bigl(
        \widetilde{Z}_N \norm{\tilde{h}_{\mathcal{A}_N, \mathcal{B}_N}^N}_{\tilde{\mu}_N}
      \bigr)
      \leq
      \alpha_N +  \frac{c_5}{N}.
    \end{align}
  \item For any $q \in [1,\infty)$ there exists a $c_6 \in (0, \infty)$ such that, for $N$ sufficiently large depending on $\beta, k_1, k_2, k_J, q$, $\prob$-a.s.\ on the event $\Otimeta(N)$,
    \begin{align}
      \label{eq:moments:harm:annealed}
      \me^{\alpha_N}(1 - c_6 N^{-1})
      \leq
      \frac{
        \mean_{\mathcal{G}}\Bigl[
          \bigl(Z_N \norm{h_{\mathcal{A}_N, \mathcal{B}_N}^N}_{\mu_N}\bigr)^q
        \Bigr]^{1/q}
      }{
        \widetilde{Z}_N 
        \norm{\tilde{h}_{\mathcal{A}_N, \mathcal{B}_N}^N}_{\tilde{\mu}_N}
      }
      \leq
      \me^{q \alpha_N} (1 + c_6 N^{-1}).
    \end{align}
  \end{enumerate}
\end{proposition}

\begin{proof}
  \textit{(i)} By using Jensen's inequality, $\cG$-measurability of $\cS_{i,N}$ and Lemma~\ref{lem:EB:Delta+HvH}(i), we find that for $2 \leq N \in \mathbb{N}$ 
  \begin{align}
    \mean_{\mathcal{G}}
    \Bigl[
      \log\bigl(Z_N \mu_N[\cS_{i,N}]\bigr)
    \Bigr]
    &\leq
    \log\mean_{\mathcal{G}}\Bigl[ Z_N \mu_N[\cS_{i,N}] \Bigr]
    \nonumber\\ 
    &=
    \log\bigl(\widetilde{Z}_N \tilde{\mu}_N[\cS_{i,N}]\bigr)
    +
    \log\biggl(
      \sum_{\sigma \in \cS_{i,N}} 
      \frac{\tilde{\mu}_N(\sigma)}{\tilde{\mu}_N[\cS_{i,N}]} 
      \mean_{\mathcal{G}}\biggl[\me^{-\beta \Delta_N(\sigma)}\biggr]
    \biggr)
    \nonumber\\
    &=
    \log\bigl(\widetilde{Z}_N \tilde{\mu}_N[\cS_{i,N}]\bigr) 
    + \alpha_N + \log\bigl(1 + O(N^{-1})\bigr).
    \label{eq:log:Zh:ub}
  \end{align}
  Likewise, 
  \begin{align}
    \mean_{\mathcal{G}}\Bigl[
      \log\bigl(Z_N \mu_N[\cS_{i,N}]\bigr)
    \Bigr]
    &=
    \log\bigl(\widetilde{Z}_N \tilde{\mu}_N[\cS_{i,N}]\bigr)
    +
    \mean_{\mathcal{G}}\biggl[\log\biggl(
        \sum_{\sigma \in \cS_{i,N}} 
        \frac{\tilde{\mu}_N(\sigma)}{\tilde{\mu}_N[\cS_{i,N}]} \,
        \me^{-\beta \Delta_N(\sigma)}
      \biggr)
    \biggr]
    \nonumber\\
    &\geq
    \log\bigl(\widetilde{Z}_N \tilde{\mu}_N[\cS_{i,N}]\bigr)
    +
    \sum_{\sigma \in \cS_{i,N}} 
    \frac{\tilde{\mu}_N(\sigma)}{\tilde{\mu}_N[\cS_{i,N}]}
    \mean_{\mathcal{G}}\bigl[-\beta \Delta_N(\sigma)\bigr]
    \nonumber\\
    &=
    \log\bigl(\widetilde{Z}_N \tilde{\mu}_N[\cS_{i,N}]\bigr).
    \label{eq:log:Zh:lb}
  \end{align}
  Moreover, since $\tilde{\mu}_N[\cS_{j,N}] / \tilde{\mu}_N[\cS_{i,N}] \leq \me^{\beta (k_J+h) N}$ for all $j \in \{1, \ldots, K\}$, we deduce from \cite[Theorem~1.7]{SS19} that
  \begin{align}
    \label{eq:tilde:harm:sum}
    \widetilde{Z}_N 
    \norm{\tilde{h}_{\mathcal{A}_N, \mathcal{B}_N}^N}_{\tilde{\mu}_N}
    =
    \widetilde{Z}_N \tilde{\mu}_N\bigl[\cS_{i,N}\bigr]
    \bigl(1 + O(\me^{-k_2N} + N \tilde{\rho}_N)\bigr).
  \end{align}
  Thus, recalling that $\tilde{\rho}_N= \me^{-k_1 N}$, the assertion in \eqref{eq:log:harm:annealed} follows from Corollary~\ref{cor:mean:harm:estimate} combined with \eqref{eq:log:Zh:ub}, \eqref{eq:log:Zh:lb} and \eqref{eq:tilde:harm:sum}.
  
  \noindent
  \textit{(ii)} For given $q \in [1,\infty)$, take $c' \in (0,\infty)$ and 
  \begin{align*}
    a_N \coloneqq k_J\sqrt{(N+1) \log 2 + q N (\log 2 + 2 \beta (k_J+h) + c')}.
  \end{align*}
  It follows from Lemma~\ref{lemma:conc:H_N} that, $\prob$-a.s.,
  \begin{align}
    \label{eq:PXic1/q}
    \prob_{\mathcal{G}}\bigl[\Xi(a_N)^c\bigr]^{\frac{1}{q}}
    \leq
    \me^{-b_N/q}
    =
    2^{-N} \me^{-2 \beta (k_J+h) N} \me^{-c' N}.
  \end{align}
  To get an upper bound for the $q$-th conditional moment given $\mathcal{G}$ of the harmonic sum, we use Minkowski's inequality, \eqref{eq:PXic1/q} and the facts that $Z_N \norm{h_{\mathcal{A}_N, \mathcal{B}_N}^N}_{\mu_N} \leq 2^N \me^{\beta (k_J+h) N}$ and $\widetilde{Z}_N \norm{\tilde{h}_{\mathcal{A}_N, \mathcal{B}_N}^N}_{\tilde{\mu}_N} \geq \me^{-\beta (k_J+h) N}$. This gives
  \begin{align}
    \mean_{\mathcal{G}}&\Bigl[
      \bigl(Z_N \norm{h_{\mathcal{A}_N, \mathcal{B}_N}^N}_{\mu_N}\bigr)^q
    \Bigr]^{\frac{1}{q}}
    \nonumber\\ 
    &\leq
    \mean_{\mathcal{G}}\Bigl[
      \bigl(Z_N \norm{h_{\mathcal{A}_N, \mathcal{B}_N}^N}_{\mu_N}\bigr)^q
      \ind_{\Xi(a_N)}
    \Bigr]^{\frac{1}{q}}
    +
    2^N \me^{\beta (k_J+h) N}
    \prob_{\mathcal{G}}\bigl[\Xi(a_N)^c\bigr]^{\frac{1}{q}}
    \nonumber\\ 
    &\leq
    \mean_{\mathcal{G}}\Bigl[
      \bigl(Z_N \norm{h_{\mathcal{A}_N, \mathcal{B}_N}^N}_{\mu_N}\bigr)^q
      \ind_{\Xi(a_N)}
    \Bigr]^{\frac{1}{q}}
    +
    \widetilde{Z}_N 
    \norm{\tilde{h}_{\mathcal{A}_N, \mathcal{B}_N}^N}_{\tilde{\mu}_N}\,
    \me^{-c' N}.
    \label{eq:moments:harm:est:ub:1}
  \end{align}
  To analyse the first term of the right-hand side of \eqref{eq:moments:harm:est:ub:1}, we apply Proposition~\ref{prop:harm:estimate} to obtain that there exists a $C \in (0, k_1 \wedge k_2)$ such that, for $N$ sufficiently large depending on $\beta, k_1, k_2, k_J, q$, 
  \begin{align}
    \label{eq:moments:harm:est:ub:2}
    \mean_{\mathcal{G}}\Bigl[
      \bigl(Z_N \norm{h_{\mathcal{A}_N, \mathcal{B}_N}^N}_{\mu_N}\bigr)^q
      \ind_{\Xi(a_N)}
    \Bigr]^{\frac{1}{q}}
    =
    \mean_{\mathcal{G}}\Bigl[
      \bigl(Z_N \mu_N\bigl[\cS_{i,N}\bigr]\bigr)^q \ind_{\Xi(a_N)}
    \Bigr]^{\frac{1}{q}} 
    \bigl(1 + O(\me^{-C N})\bigr).
  \end{align}
  Moreover, a further application of Minkowski's inequality yields that
  \begin{align}
    \mean_{\mathcal{G}}\Bigl[
      \bigl(Z_N \mu_N\bigl[\cS_{i,N}\bigr]\bigr)^q \ind_{\Xi(a_N)}
    \Bigr]^{\frac{1}{q}}
    &\leq
    \widetilde{Z}_N \sum_{\sigma \in \cS_{i,N}} \tilde{\mu}_N(\sigma)
    \mean_{\mathcal{G}}\Bigl[ \me^{-\beta q \Delta_N(\sigma)} \Bigr]^{\frac{1}{q}}
    \nonumber\\
    &=
    \widetilde{Z}_N \tilde{\mu}_N[\cS_{i,N}]\, \me^{q \alpha_N} 
    \bigl(1 + O(N^{-1})\bigr),
    \label{eq:moments:harm:est:ub:3}
  \end{align}
  where in the last step we used Lemma~\ref{lem:EB:Delta+HvH}(i) with $\beta$ replaced by $\beta q$. Thus, combining \eqref{eq:moments:harm:est:ub:1} with \eqref{eq:moments:harm:est:ub:2}, \eqref{eq:moments:harm:est:ub:3} and $\eqref{eq:tilde:harm:sum}$, we see that there exists a $c'' \in (0, \infty)$ such that for $N$ sufficiently large depending on $\beta, k_1, k_2, k_J, q$, 
  \begin{align}
    \mean_{\mathcal{G}}\Bigl[
      \bigl(Z_N \norm{h_{\mathcal{A}_N, \mathcal{B}_N}^N}_{\mu_N}\bigr)^q
    \Bigr]^{\frac{1}{q}}
    \leq
    \widetilde{Z}_N 
    \norm{\tilde{h}_{\mathcal{A}_N, \mathcal{B}_N}^N}_{\tilde{\mu}_N}\,
    \me^{q \alpha_N}
    \bigl(1 + c'' N^{-1}\bigr).
  \end{align}
  
  We close by proving a lower bound for the $q$-th conditional moment given $\mathcal{G}$ of the harmonic sum. By Jensen's inequality we get that
  \begin{align}
    \label{eq:moments:harm:est:lb:1}
    \mean_{\mathcal{G}}\Bigl[
      \bigl(Z_N \norm{h_{\mathcal{A}_N, \mathcal{B}_N}^N}_{\mu_N}\bigr)^q
    \Bigr]^{\frac{1}{q}}
    \geq
    \mean_{\mathcal{G}}\Bigl[
      Z_N \norm{h_{\mathcal{A}_N, \mathcal{B}_N}^N}_{\mu_N}
    \Bigr]
    \geq
    \mean_{\mathcal{G}}\Bigl[
      Z_N \norm{h_{\mathcal{A}_N, \mathcal{B}_N}^N}_{\mu_N} \ind_{\Xi(a_N)}
    \Bigr].
  \end{align}
  In view of Proposition~\ref{prop:harm:estimate}, together with \eqref{eq:PXic1/q} and the facts that  $ Z_N \mu_N[\cS_{i,N}]  \leq 2^N \me^{\beta (k_J+h) N}$ and $\widetilde{Z}_N \norm{\tilde{h}_{\mathcal{A}_N, \mathcal{B}_N}^N}_{\tilde{\mu}_N} \geq \me^{-\beta (k_J+h) N}$, we get that there exists a $C \in (0, k_1 \wedge k_2)$ such that, for $N$ sufficiently large depending on $\beta, k_1, k_2, k_J, q$
  \begin{align}
    \label{eq:moments:harm:est:lb:2}
    \begin{split}
      \mean_{\mathcal{G}}
      &\Bigl[
        Z_N \norm{h_{\mathcal{A}_N, \mathcal{B}_N}^N}_{\mu_N} \ind_{\Xi(a_N)}
      \Bigr]
      \\ 
      &=
      \mean_{\mathcal{G}}\Bigl[ Z_N \mu_N[\cS_{i,N}] \ind_{\Xi(a_N)} \Bigr]
      \bigl(1 + O(\me^{-C N})\bigr)
      \\ 
      &\geq
      \Bigl(
        \mean_{\mathcal{G}}\Bigl[ Z_N \mu_N[\cS_{i,N}] \Bigr] 
        - \widetilde{Z}_N 
        \norm{\tilde{h}_{\mathcal{A}_N, \mathcal{B}_N}^N}_{\tilde{\mu}_N}\,
        \me^{-c' N}
      \Bigr)
      \bigl(1 + O(\me^{-C N})\bigr).
    \end{split}
  \end{align}
  Since, by Lemma~\ref{lem:EB:Delta+HvH}(i),
  \begin{align}
    \label{eq:moments:harm:est:lb:3}
    \mean_{\mathcal{G}}\Bigl[ Z_N \mu_N[\cS_{i,N}] \Bigr]
    =
    \widetilde{Z}_N \tilde{\mu}_N[\cS_{i,N}]\, \me^{\alpha_N} 
    \bigl(1 + O(N^{-1})\bigr),
  \end{align}
  we conclude from \eqref{eq:moments:harm:est:lb:1} combined with \eqref{eq:moments:harm:est:lb:2}, \eqref{eq:moments:harm:est:lb:3} and \eqref{eq:tilde:harm:sum} that there exists a $c''' \in (0, \infty)$ such that, for $N$ sufficiently large depending on $\beta, k_1, k_2, k_J, q$
  \begin{align}
    \mean_{\mathcal{G}}\Bigl[
      \bigl(Z_N \norm{h_{\mathcal{A}_N, \mathcal{B}_N}^N}_{\mu_N}\bigr)^q
    \Bigr]^{\frac{1}{q}}
    \geq
    \widetilde{Z}_N 
    \norm{\tilde{h}_{\mathcal{A}_N, \mathcal{B}_N}^N}_{\tilde{\mu}_N}\,
    \me^{\alpha_N}
    \bigl(1 + c''' N^{-1}\bigr).
  \end{align}
  By setting $c_6 \coloneqq c'' \vee c'''$, we get the assertion.
\end{proof}


\section{Estimates on mean hitting times of metastable sets} \label{sec:mean:hitting:times}

Before proving Theorem~\ref{thm_main}, we state two immediate corollaries of the propositions proved in Sections~\ref{sec:capa} and \ref{sec:harm}.

\begin{corollary}
  \label{cor:conc:log:MHT}
  There exist $C \in (0, k_1 \wedge k_2)$  and $c_4 \in (0, \infty)$ such that, for all $N$ sufficiently large depending on $\beta, k_1, k_2, k_J$, and all $t \in \mathbb{N}_0$, $\prob$-a.s.,
  \begin{align}
    \label{eq:conc:log:MHT}
    \begin{split}
      \prob_{\cG}\biggl[
        \abs{
          \log
          \bigl(
            \Mean_{\nu_{\mathcal{A}_N,\mathcal{B}_N}}^N\bigl[\tau_{\mathcal{B}_N}^N\bigr]
          \bigr)
          -
          \mean_{\cG}\Bigl[
            \log
            \bigl( 
              \Mean_{\nu_{\mathcal{A}_N,\mathcal{B}_N}}^N\bigl[
                \tau_{\mathcal{B}_N}^N
              \bigr]
            \bigr)
          \Bigr]
        } > t,
        \Otimeta(N)
      \biggr]
      \\ 
      \leq
      \ind_{\Otimeta(N)}
      \biggl[
        2 \biggl(\me^{-\bigl(\frac{t-c_N}{2\beta k_J}\bigr)^2 } + \me^{-\bigl(\frac{t}{2\beta k_J}\bigr)^2}\biggr) +  \me^{-k_1 N}
      \biggr],
    \end{split}
  \end{align}
  where $c_N \coloneqq c_4 \,\me^{-C N}$.
\end{corollary}

\begin{proof}
  In view of \eqref{eq:PotTheor}, we have that
  \begin{align}
    \begin{split}
      &\prob_{\cG}
      \biggl[
        \abs{
          \log
          \bigl(
            \Mean_{\nu_{\mathcal{A}_N,\mathcal{B}_N}}^N\bigl[
              \tau_{\mathcal{B}_N}^N
            \bigr]
          \bigr)
          -
          \mean_{\cG}\Bigl[
            \log
            \bigl( 
              \Mean_{\nu_{\mathcal{A}_N,\mathcal{B}_N}}^N\bigl[
                \tau_{\mathcal{B}_N}^N
              \bigr]
            \bigr)
          \Bigr]
        } > t,
        \Otimeta(N)
      \biggr]
      \\ 
      &\leq
      \prob_{\cG}\biggl[
        \abs{
          \log \bigl(Z_N \norm{h_{\mathcal{A}_N, \mathcal{B}_N}^N}_{\mu_N}\bigr)
          -
          \mean_{\cG}\Bigl[
            \log 
            \bigl(Z_N \norm{h_{\mathcal{A}_N, \mathcal{B}_N}^N}_{\mu_N}\bigr)
          \Bigr]
        } > \frac{t}{2}, 
        \Otimeta(N)
      \biggr]
      \\
      &\,\,+
      \prob_{\cG}\biggl[
        \abs{
          \log \bigl(Z_N \capa_N(\mathcal{A}_N, \mathcal{B}_N)\bigr)
          -
          \mean_{\cG}\Bigl[
            \log \bigl(Z_N \capa_N(\mathcal{A}_N, \mathcal{B}_N)\bigr)
          \Bigr]
        } > \frac{t}{2},
        \Otimeta(N)
      \biggr].
    \end{split}
  \end{align}
  Thus, the assertion follows immediately from Proposition~\ref{pro:conc_logcap} and Proposition~\ref{pro:conc_Zh}.
\end{proof}

\begin{corollary}
  \label{cor:annealed:log:MHT:comparison}
  There exists a $c_7 \in (0, \infty)$ such that, for $N$ sufficiently large depending on $\beta, k_1, k_2, k_J$,  $\prob$-a.s.\ on the event $\Otimeta(N)$,
  \begin{align}
    \label{eq:annealed:log:MHT:comparison}
    -\alpha_N - \frac{c_7}{\sqrt{N}}
    \leq
    \mean_{\mathcal{G}}\Bigl[
      \log
      \Mean_{\nu_{\mathcal{A}_N, \mathcal{B}_N}}^N\bigl[
        \tau_{\mathcal{B}_N}^N
      \bigr]
    \Bigr]
    -
    \log
    \tMean_{\tilde{\nu}_{\mathcal{A}_N, \mathcal{B}_N}}^N\bigl[
      \tilde{\tau}_{\mathcal{B}_N}^N
    \bigr]
    \leq
    2\alpha_N +\frac{c_7}{\sqrt{N}},
  \end{align}
  where $\alpha_N$ is defined in \eqref{eq:defAlphaN}.
\end{corollary}

\begin{proof}
  In view of \eqref{eq:PotTheor}, the assertion in \eqref{eq:annealed:log:MHT:comparison} follows immediately from Proposition~\ref{prop:harm:annealed}(i) and Proposition~\ref{prop:EB:cap:estimates}(i).
\end{proof}

\begin{proof}[Proof of Theorem~{\upshape\ref{thm_main}}]
  \textit{(i)} 
  Recall once again the $\cG$-measurability of $\Otimeta(N)$ (see Remark~\ref{rem:Gmeas}), and note that an application of Corollary~\ref{cor:annealed:log:MHT:comparison} yields that there exists a $c_7 \in (0, \infty)$ such that, for $N$ sufficiently large depending on $\beta, k_1, k_2, k_J$, $\prob$-a.s.,
  \begin{align}
    &\prob_{\mathcal{G}}
    \Biggl[
      \me^{-t - \alpha_N} 
      \leq
      \frac{
        \Mean_{\nu_{\mathcal{A}_N, \mathcal{B}_N}}^N\bigl[
          \tau_{\mathcal{B}_N}^N
        \bigr]
      }{
        \tMean_{\tilde{\nu}_{\mathcal{A}_N, \mathcal{B}_N}}^N\bigl[
          \tilde{\tau}_{\mathcal{B}_N}^N
        \bigr]
      }
      \leq
      \me^{+t + 2\alpha_N}
    \Biggr] 
    \ind_{\Otimeta(N)}
    \nonumber\\ 
    &\overset{\mspace{-12mu}\eqref{eq:annealed:log:MHT:comparison}}{\geq}\,
    \prob_{\mathcal{G}}\Biggl[
      \me^{-(t - \frac{c_7}{\sqrt{N}})} 
      \leq
      \frac{
        \Mean_{\nu_{\mathcal{A}_N, \mathcal{B}_N}}^N\bigl[
          \tau_{\mathcal{B}_N}^N
        \bigr]
      }{
        \exp\Bigl(
          \mean_{\mathcal{G}}\Bigl[
            \log
            \Mean_{\nu_{\mathcal{A}_N, \mathcal{B}_N}}^N\bigl[
              \tau_{\mathcal{B}_N}^N
            \bigr]
          \Bigr]
        \Bigr)
      }
      \leq
      \me^{+(t - \frac{c_7}{\sqrt{N}})}, \Otimeta(N)
    \Biggr]
    \nonumber\\ 
    &\geq
    \ind_{\Otimeta(N)} 
    \nonumber\\
    &\quad -
    \prob_{\mathcal{G}}\biggl[
      \abs{
        \log
        \bigl(
          \Mean_{\nu_{\mathcal{A}_N,\mathcal{B}_N}}^N\bigl[\tau_{\mathcal{B}_N}^N\bigr]
        \bigr)
        -
        \mean_{\mathcal{G}}\Bigl[
          \log
          \bigl( 
            \Mean_{\nu_{\mathcal{A}_N,\mathcal{B}_N}}^N\bigl[
              \tau_{\mathcal{B}_N}^N
            \bigr]
          \bigr)
        \Bigr]
      } > t - \frac{c_7}{\sqrt{N}}
    \biggr]
    \ind_{\Otimeta(N)}.
  \end{align}
  Thus, from Corollary~\ref{cor:conc:log:MHT} it follows that there exists a $C \in (0, k_1 \wedge k_2)$ such that, for $N$ sufficiently large depending on $\beta, k_1, k_2, k_J$, $\prob$-a.s.,
  \begin{align}
    \prob_{\cG}&\Biggl[
      \me^{-t - \alpha_N} 
      \leq
      \frac{
        \Mean_{\nu_{\mathcal{A}_N, \mathcal{B}_N}}^N\bigl[
          \tau_{\mathcal{B}_N}^N
        \bigr]
      }{
        \tMean_{\tilde{\nu}_{\mathcal{A}_N, \mathcal{B}_N}}^N\bigl[
          \tilde{\tau}_{\mathcal{B}_N}^N
        \bigr]
      }
      \leq
      \me^{+t + 2\alpha_N}
    \Biggr]
    \ind_{\Otimeta(N)}
    \nonumber\\ 
    &\geq 
    \biggl(
      1 -
      2 \Bigl(\me^{-\bigl(\frac{t - c_7 N^{-1/2}-c_N}{2\beta k_J}\bigr)^2 } + \me^{-\bigl(\frac{ t - c_7 N^{-1/2}}{2\beta k_J}\bigr)^2}\Bigr) -  \me^{-k_1 N}\biggr)
    \ind_{\Otimeta(N)}.
    \label{eq:lastStep}
  \end{align} 
  Assumption~\ref{ass:meta} implies that $\lim_{N \to \infty} \ind_{\Otimeta(N)}$ exists and is $\prob$-a.s.\ equal to $1$. Since trivially $1 \geq \ind_{\Otimeta(N)}$, taking the limit $N \to \infty$ of \eqref{eq:lastStep} yields \eqref{eq:main:conc}.
  
  \noindent
  \textit{(ii)} Fix $q \in [1, \infty)$.  In view of \eqref{eq:PotTheor}, an application of the Cauchy-Schwarz inequality yields that, $\prob$-a.s.,
  \begin{align}
    \mean_{\mathcal{G}}&\Bigl[
      \Mean_{\nu_{\mathcal{A}_N, \mathcal{B}_N}}^N\bigl[
        \tau_{\mathcal{B}_N}^N
      \bigr]^q
    \Bigr]^{\frac{1}{q}}
    \nonumber\\
    &\leq
    \mean_{\mathcal{G}}\Bigl[
      \bigl(
        Z_N \norm{h_{\mathcal{A}_N, \mathcal{B}_N}^N}_{\mu_N}
      \bigr)^{2q}
    \Bigr]^{\frac{1}{2q}}
    \mean_{\mathcal{G}}\Bigl[
      \bigl(
        Z_N \capa_N\bigl(\mathcal{A}_N, \mathcal{B}_N\bigr)
      \bigr)^{-2q}
    \Bigr]^{\frac{1}{2q}}.
  \end{align}
  Hence, by Proposition~\ref{prop:EB:cap:estimates}(ii) and \ref{prop:harm:annealed}(ii), there exists a $c \in (0, \infty)$ such that for $N$ sufficiently large depending on $\beta, k_1, k_2, k_J, q$, $\prob$-a.s.\ on the event $\Otimeta(N)$,
  \begin{align}
    \mean_{\mathcal{G}}\Bigl[
      \Mean_{\nu_{\mathcal{A}_N, \mathcal{B}_N}}^N\bigl[
        \tau_{\mathcal{B}_N}^N
      \bigr]^q
    \Bigr]^{\frac{1}{q}}
    \leq
    \tMean_{\tilde{\nu}_{\mathcal{A}_N, \mathcal{B}_N}}^N\bigl[
      \tilde{\tau}_{\mathcal{B}_N}^N
    \bigr]\,
    \me^{4q \alpha_N} \Bigl(1 + \frac{c}{\sqrt{N}}\Bigr).
  \end{align}
  On the other hand, by using Jensen's inequality and Corollary~\ref{cor:annealed:log:MHT:comparison}, we find that there exists a $c_7 \in (0, \infty)$ such that, for any $q \in [1, \infty)$ and $N$ sufficiently large depending on $\beta, k_1, k_2, k_J,q$, $\prob$-a.s.\ on the event $\Otimeta(N)$,
  \begin{align}
    \mean_{\mathcal{G}}\Bigl[
      \Mean_{\nu_{\mathcal{A}_N, \mathcal{B}_N}}^N\bigl[
        \tau_{\mathcal{B}_N}^N
      \bigr]^q
    \Bigr]^{\frac{1}{q}}
    &\geq
    \exp\Bigl(
      \mean_{\mathcal{G}}\Bigl[
        \log
        \Mean_{\nu_{\mathcal{A}_N, \mathcal{B}_N}}^N\bigl[
          \tau_{\mathcal{B}_N}^N
        \bigr]
      \Bigr]
    \Bigr)
    \nonumber\\ 
    &\geq
    \tMean_{\tilde{\nu}_{\mathcal{A}_N, \mathcal{B}_N}}^N\bigl[
      \tilde{\tau}_{\mathcal{B}_N}^N
    \bigr]\,
    \me^{-\alpha_N} \Bigl(1 - \frac{c_7}{\sqrt{N}}\Bigr).
  \end{align}
  Therefore, by letting $c_1= c \vee c_7$,  for $N$ sufficiently large depending on $\beta, k_1, k_2, k_J,q$, $\prob$-a.s.\ on the event $\Otimeta(N)$, 
  \begin{align}
    \me^{-\alpha_N}\bigl(1 - \tfrac{c_1}{\sqrt{N}}\bigr)
    \leq
    \frac{
      \mean_{\cG}\Bigl[
        \Mean_{\nu_{\mathcal{A}_N, \mathcal{B}_N}}^N\bigl[\tau_{\mathcal{B}_N}^N\bigr]^q
      \Bigr]^{1/q}
    }
    {
      \tMean_{\tilde{\nu}_{\mathcal{A}_N, \mathcal{B}_N}}^N\bigl[\tilde{\tau}_{\mathcal{B}_N}^N\bigr]
    }
    \leq
    \me^{4q \alpha_N}\bigl(1 + \tfrac{c_1}{\sqrt{N}}\bigr).
  \end{align}
  Thus, the set $\Omega_{q,c_1}(N)$ defined in \eqref{eq:def_Oqc} contains $\Otimeta(N)$. Therefore, using Assumption~\ref{ass:meta}, and monotonicity of probability
  \begin{align}
    1= 
    \prob\biggl[
      \liminf_{N \to \infty} \Otimeta(N)
    \biggr]
    \leq \prob\biggl[
      \liminf_{N \to \infty} \Omega_{q,c_1}(N)
    \biggr]
  \end{align}
  which concludes the proof of \eqref{eq:main:moments}.
\end{proof}  


\appendix

\section{Concentration inequality} \label{app:McD}

We present a concentration inequality for functionals of conditionally independent random variables that is a slight extension of the classical McDiarmid concentration inequality for functionals of independent random variables satisfying a bounded difference estimate (cf.\ \cite[Theorem~6.2]{BLM13}, \cite[Section~2.4.1]{DZ10}).

\begin{proposition}
  \label{prop:McDiarmid}
  Let $(\Omega, \mathcal{F}, \prob)$ be a probability space, $\mathcal{G} \subset \mathcal{F}$ a sub-$\sigma$-algebra of $\mathcal{F}$, $1\leq n \in \mathbb{N}$ and $\mathcal{X}$ a Polish space. Consider a vector $X = (X_1, \ldots, X_n)$ of $\mathcal{X}$-valued random variables on $(\Omega, \mathcal{F}, \prob)$ that are conditionally independent given $\mathcal{G}$, and let $f_n : \mathcal{X}^n \to \mathbb{R}$ be a measurable function.  Suppose that, for any $i \in \{1, \ldots, n\}$, 
  \begin{align}
    \label{eq:bounded:difference}
    \bigl\vert 
      f_n(X_1, \ldots, X_n) 
      - f_n(X_1, \ldots, X_{i-1}, X_i', X_{i+1}, \ldots, X_n)
    \bigr\vert 
    \leq
    c_i \in [0, \infty) \qquad \prob\text{-a.s.},
  \end{align}
  where $(X_1', \ldots, X_n')$ is a conditionally independent copy of $(X_1, \ldots, X_n)$ given $\mathcal{G}$.  Then, $\prob$-a.s.\ for all $t \in \mathbb{N}_0$,
  \begin{align}
    \left.
      \begin{array}{c}
        \prob\bigl[
          f_n(X) - \mean[f_n(X) \mid  \mathcal{G}] > +t \mid \mathcal{G}
        \bigr]
        \\[1ex]
        \prob\bigl[
          f_n(X) - \mean[f_n(X) \mid \mathcal{G}] < -t \mid \mathcal{G}
        \bigr]
      \end{array}
    \right\}
    \leq
    \me^{-t^2/(2v)},
  \end{align}
  where $v \coloneqq \frac{1}{4} \sum_{i=1}^n c_i^2$.
\end{proposition}  

\begin{proof}
  Since there exists a regular conditional probability for $X$ (see e.g.\ \cite[p.~217]{CT97}), the proof follows the line of proof of the non-conditional McDiarmid concentration inequality.
\end{proof}



\end{document}